\documentclass[11pt]{article}
\usepackage[margin=1in]{geometry}
\usepackage{amsmath,amssymb,amsthm}
\usepackage{commath}
\usepackage{hyperref}
\usepackage{algorithm}
\usepackage{color}
\usepackage{graphicx}
\usepackage{algpseudocode}

\DeclareMathOperator*{\argmin}{arg\,min}

\newcommand\Rbb{\mathbb{R}}

\newcommand\irm{\mathrm{i}}

\newcommand\Fcal{\mathcal{F}}

\newcommand\Dhat{\hat{D}}

\title{Joint image edge reconstruction and its application in multi-contrast MRI}
\author{Yunmei Chen
\footnote{Department of Mathematics, University of Florida, Gainesville, FL 32611, USA. Email: \texttt{yun@math.ufl.edu}}
\and
Ruogu Fang
\footnote{Department of Biomedical Engineering, University of Florida, Gainesville, FL 32611, USA. Email: \texttt{ruogu.fang@bme.ufl.edu}} 
\and
Xiaojing Ye
\footnote{Corresponding author. Department of Mathematics \& Statistics, Georgia State University, Atlanta, GA 30303, USA.
Email: \texttt{xye@gsu.edu}}
}
\date{}

\begin{document}
\maketitle

\begin{abstract}	
We propose a new joint image reconstruction method by recovering edge directly from observed data. More specifically, we reformulate joint image reconstruction with vectorial total-variation regularization as an $l_1$ minimization problem of the Jacobian of the underlying multi-modality or multi-contrast images. Derivation of data fidelity for Jacobian and transformation of noise distribution are also detailed. The new minimization problem yields an optimal $O(1/k^2)$ convergence rate, where $k$ is the iteration number, and the per-iteration cost is low thanks to the close-form matrix-valued shrinkage. We conducted numerical tests on a number multi-contrast magnetic resonance image (MRI) datasets, which show that the proposed method significantly improves reconstruction efficiency and accuracy compared to the state-of-the-arts.
\bigskip

\end{abstract}

\section{Introduction}
\label{sec:intro}
The advances of medical imaging technology have allowed simultaneous
data acquisition for multi-modality images such as PET-CT, PET-MRI and multi-contrast images
such as T1/T2/PD MRI. 
Multi-modality/contrast imaging integrates two or more imaging modalities (or contrasts)
into one system to produce more comprehensive observations of the subject.
Such technology combines the strengths of different imaging modalities/contrasts 
in clinical diagnostic imaging, and hence can be much more precise and effective 
than conventional imaging.
The images from different modalities or contrasts complement each other
and generate high spatial resolution and better tissue contrast.
However, effectively utilizing sharable information from different modalities and reconstructing
multi-modality images remain as a challenging task especially when only limited data
are available.
Therefore, the main approach of joint image reconstruction is
to incorporate similarities, such as anatomical structure, between all modalities/contrasts 
into the process to improve reconstruction accuracy \cite{Bilgic:2011a,Ehrhardt:2016b,Huang:2014a}. 

In this paper, we consider a new approach by jointly reconstructing edges of images, 
rather than the images themselves, directly from multi-modality/contrast imaging data. 
The final images can be obtained very easily given the reconstructed edges.
We show that this new
approach results in an $l_1$-type minimization that can be effectively solved by
accelerated prox-gradient methods with an optimal $O(1/k^2)$ convergence rate,
where $k$ is the iteration number. Moreover, the subproblems reduce complex
vectorial/joint total variation regularizations to simple matrix-valued shrinkages,
which often have cheap closed-form solutions. This is in sharp contrast to 
primal-dual based image reconstruction algorithms, such as primal-dual hybrid gradient
method (PDHG) and alternating direction
method of multipliers (ADMM), which only yield $O(1/k)$
convergence rate. Therefore, the proposed method can produce high quality
reconstructions with much improved efficiency over the state-of-the-arts methods
in multi-modality/contrast image reconstruction problems.

The contributions of this paper can be summarized as follows.
(a) We develop a novel two-step joint image reconstruction method
that transforms the vectorial TV regularized minimization of image into
an $l_1$ minimization of Jacobian. 
This enables a numerical scheme with optimal $O(1/k^2)$ convergence rate
by employing accelerated gradient descent method
(Section \ref{subsec:edgerecon}).
(b) In the resulting $l_1$ minimization problems, the main subproblem involves
matrix norm (instead of vectorial TV) and can 
be solved easily using generalized shrinkage for matrices.
We provide close-form solutions for several cases originated from commonly used vectorial TVs
(Section \ref{subsec:shrinkage} and Appendix \ref{apd:shrinkage}).
(c) We analyze the noise distribution after transformation, and incorporate it into the data fidelity of
Jacobian in the algorithm which significantly improves reconstruction accuracy and efficiency 
(Section \ref{subsec:data_reform}).
(d) We conduct a series to numerical tests on several multi-contrast MRI datasets
and show the very promising performance of the proposed methods (Section \ref{sec:results}).
Although we only focus on multi-contrast MRI reconstruction where all data is acquired in Fourier space, 
it is worth noting that the proposed method can be readily extended to other cases, such as those
with Radon data, 
hence is applicable to reconstructions involving other types of imaging modalities.

The remainder of the paper is organized as follows. We first provide an overview
of recent literatures in joint image reconstruction in Section \ref{sec:related}.
Then we present the proposed method and address a number of details
in Section \ref{sec:proposed}. In Section \ref{sec:results}, we conduct
a number of numerical tests on a variety of multi-contrast MR image
reconstruction datasets. Section \ref{sec:conclusion} concludes our findings.

\section{Related Work}
\label{sec:related}
There have been a significant amount effort 
devoted to develop models and algorithms that 
can effectively take the anatomy structure similarities across modalities/contrasts 
into account during joint reconstructions. 
As widely accepted, these structure similarities can be exploited using the locations and directions of the edges 
\cite{Ehrhardt:2015b,Ehrhardt:2016a,Ehrhardt:2015a} characterized by the magnitude 
and direction of the gradient of an image. The active researches on multi-modal/contrast image reconstructions 
have focused  mainly on how to effectively utilize the complimentary information on these structure similarities 
to improve the accuracy and robustness of the reconstructions.

Inspired by the success of total variation (TV) based image reconstructions
for scalar-valued images, many algorithms for joint reconstruction of 
multi-modal images extend the classical TV to vectorial TV 
for joint multi-modal image reconstructions. This extension aims at 
capturing the sharable edge information and performing smoothing along the common 
edges across the modalities. There have been several different ways to define the 
TV regularization for multi-modal images (represented by vector-valued functions).
For instance, $\sum_{j=1}^m TV(u_j)$, where $u_j$ is the (scalar-valued) image
of modality/channel/contrast $j$ in an $m$-modality imaging problem, 
is proposed in \cite{Blomgren:1998a}. 
In \cite{Duran:2016a}, a number of variations along this direction, 
called collaborative TV (CTV), are studied and summarized comprehensively. 
A specific TV takes a particular form to integrate 
partial derivatives of the image across the modalities. For example,
$\sum_{j=1}^m TV(u_j)$ takes $l_2$ norm of 
image gradient at every point (pixel) for each modality,
then $l_1$ norm over all pixels, and finally 
$l_1$ norm (direct sum) across all modalities. 
Another commonly used CTV variation, called joint TV (JTV), is formulated as
$\int_\Omega (\sum_{j=1}^m \sum_{l=1}^d |\partial_l u_j|^2)^{1/2} \dif x$
and has been successfully applied to color image reconstruction \cite{Bresson:2008a,Sapiro:1996b},
multi-contrast MR image reconstruction \cite{Huang:2014a, Majumdar:2011a},
joint PET-MRI reconstruction \cite{Ehrhardt:2015a}, and in geophysics \cite{Haber:2013a}.
As one can see, JTV is taking $l_2$ norm in terms of gradients, 
$l_2$ norm across modalities, and then finally $l_1$ for pixels.
In general, the gradient of a vector-valued image consisting multiple modalities/contrasts
is a tensor at each pixel, and the TV takes specific 
combination of norms regarding 
partial derivatives (gradients), modalities (channels, or contrasts), and pixels, respectively.

Another joint regularization approach is to
extend anisotropic TV regularization from scalar to vector-valued case. 
In this approach, anisotropic TV is employed to
replace isotropic TV by an anisotropic term that 
incorporates directional structures exhibited by either the original data or the underlying image.
Anisotropic TV been applied in standard single-modal image reconstruction with successes in
\cite{Estellers:2015a,Grasmair:2010a, Lenzen:2015a}, where structure tensor is used to
provide information of both size and orientation of image gradients. 
In particular, the model proposed in \cite{Grasmair:2010a} employs an anisotropic TV regularization
of form $\int_\Omega (\dif u^T A(u)\dif u)^{1/2} \dif x$
for image reconstruction. Here, $\dif u=(\partial_1 u,\partial_2 u)^T$ is 
the gradient of the single-modal, scalar-valued 2D image
$u$, the diffusion tensor $A(u)$ is determined by the 
eigenvalues $(\lambda_1, \lambda_2)$ and corresponding eigenvectors $(v_1,v_2)$ of the 
structure tensor $J_\rho =k_\rho \dif u \dif u^T$.
The model developed in \cite{Estellers:2015a} alternately solves two minimizations:
one that estimates structure tensor using the image from previous iteration, 
and the other one improves the image with an adaptive regularizer defined from this tensor. 
The idea of anisotropic TV has been extended to vector-valued images with the anisotropy 
adopted from the gradients of multi-modal/contrast images, which incorporates similarities 
of directional structures in joint reconstruction. 
The method developed in \cite{Ehrhardt:2016b} projects the gradient in the total variation 
functional onto a predefined vector field given by the other contrast for joint reconstruction 
of multi-contrast MR images. In \cite{Ehrhardt:2016b}, a directional TV (DTV) is proposed with
the form $DTV_{u_2}(u_1)= \|P_{\zeta} (\dif u_1)\|_1$, where 
$P_{\zeta} (\dif u_1)$ is the residue of projection of 
$\dif u_1$ to $\zeta$,
and $u_1$ and $u_2$ are scalar-valued functions
each representing the image of one modality and $\zeta(x)=\dif u_2(x)/\|\dif u_2(x)\|$
for all $x\in \Omega$. In other words, this is the anisotropic diffusion 
$\|(I- \zeta\zeta^T) \dif u_1\|_1$ by using the structure tensor $\zeta \zeta^T$ of
given reference image $u_2$. 

With geometric interpretation of gradient tensor above, 
authors in \cite{Di-Zenzo:1986a} suggest to consider a 
vector-valued image as a parameterized 2-dimensional 
Riemann manifold with metric 
$G=D u^T D u$, where $u$ is a vector-valued image and $Du=[\partial_j u_k]_{j,k}$
is the $2\times 2$ Jacobian of $u$. Then the eigenvector corresponding 
to the larger eigenvalue gives the direction of the vectorial edge. Based on this framework, 
several forms of vectorial TV (VTV) have been developed. In \cite{Sapiro:1996a}, a family of VTV
formulations are suggested as the integral of $f(\lambda_+,\lambda_-)$ over the manifold, 
where $\lambda_+$ and $\lambda_-$ denote the larger and smaller eigenvalues of the metric 
$G$, respectively, 
and $f$ is a suitable scalar-valued function. A special choice of 
$f(\lambda_+, \lambda_-) = \sqrt{\lambda_+ + \lambda_-}$, i.e., the Frobenius norm of the Jacobian 
$Du$, reduce to JTV mentioned above. For another special choice of 
$f(\lambda_+,\lambda_-) = \sqrt{\lambda_+}$ , there is $VTV(u)=\int_\Omega \sigma_1(Du) \dif x$, 
where $\sigma_1(Du)$ is the largest singular value of the Jacobian $Du$ of $u$. 
This can be computed by using the dual formulation 
$\sup_{(\xi,\eta) \in K} \sum_{k=1}^d \int_\Omega u_k \mathrm{div}(\eta_k\xi)\dif x$  
with $K=C_0^1(\Omega; S^m \times S^d)$ where $S^d$ is the standard simplex in $\Rbb^d$.

Besides joint TV or tensor based regularization,
the parallelism of level sets across multi-modality/contrast images are also proposed
as joint image regularization in \cite{Ehrhardt:2014b,Ehrhardt:2015a, Ehrhardt:2014a,Haber:2013a}.
The main idea is to exploit the structural 
similarity of two images $u_1$ and $u_2$ measured by the parallelism of their gradients 
$\dif u_1$ and $\dif u_2$
at each point using
$\delta(\dif u_1,\dif u_2):= f (g(|\dif u_1||\dif u_2|)-g(\langle \dif u_1,\dif u_2 \rangle ))$, for some functions $f$ and $g$.
Then the regularization in joint reconstruction takes form
$\int_\Omega \delta(\dif u_1, \dif u_2) \dif x$
Several different choices of $f$ and $g$ have been studied in these works. 
For instance, in \cite{Ehrhardt:2014a}, $f$ and $g$ are taken as identities, or 
$f(s)=\sqrt{1+s}$ and $g(s)=s^2$. In \cite{Haber:2013a}, $f$ is the identity and 
$g(s)=s^2$. 
In \cite{Ehrhardt:2015a} the side information
on the level set, namely, the location and direction of 
the gradients of a reference image, is available to assist the reconstruction.

Another joint reconstruction approach different from aforementioned methods is 
to recover gradient information for each of the underlying multi-modal images 
from the measurements in Fourier domain, then use the recovered gradients
to reconstruct images. This approach is motivated by the idea of Bayesian compressed sensing 
and applied to joint multi-contrast MR image reconstruction in \cite{Bilgic:2011a}. 
In \cite{Bilgic:2011a}, gradients of the multi-contrasts images are reconstructed 
jointly from their measurements in Fourier space under a hierarchical Bayesian framework, 
where the joint sparsity on gradients across multi-contrast MRI is exploited by sharing 
the hyper-parameter in the their maximum a posteriori (MAP) estimations. 
Their experiments show the advantage of using joint sparsity on gradients over 
conventional sparsity. However, their method requires extensive computational cost.
A two-step gradient reconstruction of MR images is also proposed in \cite{Patel:2012a}, however, only
single-modality/contrast image is considered. 
In \cite{Patel:2012a}, the authors showed that this two-step gradient reconstruction approach allows to 
reconstruct image with fewer number of measurements than required by standard TV minimization method.

\section{Proposed Method}
\label{sec:proposed}
In this section, we propose a new joint image reconstruction method that 
first restores image edges (gradients/Jacobian) from observed data,
and then assembles the final image using these edges.
Without loss of generality, we assume all images are $2$-dimensional, i.e., 
the image domain $\Omega\subset \Rbb^2$ (all derivations below can be 
easily extended to higher dimensional images). 
For simplicity, we further assume $\Omega$ is rectangular. 
For single channel/modality/contrast case, we use function $u:\Omega \to \Rbb$ 
represents the image, such that $u(x)$ stands for the intensity of image at $x\in\Omega$.
In multi-modality case, $u:\Omega \to \Rbb^m$ where $m$ is the number of modalities.

It is also convenient to treat an image $u$ as a matrix in discretized setting which we mainly work on
in practice, and further
stack the columns of $u$ to form a single column vector in $\Rbb^n$, where $n$ is the total
number of pixels in the image.
For multi-modality case, we have $u_j\in\Rbb^n$ to represent the (discretized) image of
modality $j$ for $j=1,\dots,m$.

\subsection{Vectorial total-variation regularization}\label{subsec:VTV}
Standard total-variation (TV) regularized image reconstruction can be 
formulated as a minimization problem as follows:
\begin{equation}\label{eq:TVbased}
\min_u \alpha TV(u) + h(u)
\end{equation}
where $h$ represents the data fidelity function, e.g., $h(u)=\frac{1}{2}\|Au-b\|^2$,
for some given data sensing matrix $A$ and observed partial/noisy/blurry data $b$.
By solving \eqref{eq:TVbased}, we obtain a solution $u$ which minimizes
the sum of TV regularization term and data fidelity term with a weighting
parameter $\alpha>0$ that balances the two terms.
It is shown that TV regularization can effectively recover high quality images
with well preserved object boundaries from limited and/or noisy data.

In joint multi-modality image reconstruction, the edges of images from different modalities
are highly correlated. To take such factor into consideration, the standard 
TV regularized image reconstruction \eqref{eq:TVbased} can be simply replaced by vectorial TV (VTV) 
regularized counterpart:
\begin{equation}\label{eq:VTVbased}
\min_u \alpha VTV(u) + h(u)
\end{equation}
where $VTV(u)$ is the vectorial TV of $u$. In the case that
$u$ is continuously differentiable,
VTV is a direct extension of standard TV as an ``$l_1$ of gradient'' as
\begin{equation}\label{eq:VTV}
VTV(u)=\int_\Omega \|Du(x)\|_\star \dif x
\end{equation}
where $Du(x)\in\Rbb^{2\times m}$ is the Jacobian matrix at point $x$,
and $\star$ indicates some specific matrix norm. For example,
let $Q=[q_{ij}]$ be an $d$-by-$m$ matrix, then we may use one of the following
matrix norms as $\|\cdot\|_\star$:
\begin{itemize}
\item Frobenius norm: $\|Q\|_F=(\sum_{i,j}|q_{ij}|^2)^{1/2}$.
This essentially treats $Q$ as an $(dm)$-vector.
\item Induced 2-norm: $\|Q\|_2=\sigma_1$ where $\sigma_1$ is the largest singular value of $Q$.
This norm is advocated in \cite{Goldluecke:2010a} with a geometric interpretation when used in VTV.
\item Nuclear norm: $\|Q\|_*=\sum_{i=1}^{\min\{d,m\}}\sigma_i$ where $\sigma_1\geq \sigma_2\geq\dots\geq0$ 
are singular values of $Q$. This is a convex relaxation of matrix rank.
\end{itemize}
Obviously, there are a number of other choices for VTV due to the many
variations of matrix norms. However, in this paper, we only focus on these three norms
as they are the mostly used ones in VTV regularized image reconstructions.

It is also worth noting that, in general, the VTV norm with matrix norm $\|\cdot\|_\star$
is defined for any function $u\in L_\text{loc}^1(\Rbb^2; \Rbb^{m})$ (not necessarily differentiable)
as
\begin{equation}\label{eq:VTVdef}
VTV(u) = \sup \cbr[3]{ \int_\Omega u(x) \mathrm{div}(\xi(x)) \dif x: \|\xi(x)\|_\bullet\leq 1,\ \xi(x)\in\Rbb^{2\times m},\forall x\in \Omega}
\end{equation}
where $\|\cdot\|_\bullet$ is the dual norm of $\|\cdot\|_\star$. 
Although we would not always make use of the original VTV definition \eqref{eq:VTVdef}
in discrete setting, we show that they can help to derive closed-form soft-shrinkage with
respect to the corresponding matrix $\star$-norm as in Appendix \ref{apd:shrinkage}.

\subsection{Joint edge reconstruction}
\label{subsec:edgerecon}

The \textbf{main idea} of this paper is to \textbf{reconstruct edges (gradients/Jacobian) of multi-modality
images jointly}, and then assemble the final image from these edges.
To that end, we let $v$ denote the ``gradients'' (or Jacobian) of $u$. 
If $u$ is differentiable, then $v(x)=(v_{kj}(x)):=Du(x) \in\Rbb^{2\times m}$
is the Jacobian matrix of $u$ at $x\in \Omega$. 
For example, assuming there are three modalities and 
$u(x)=(u_1(x),u_2(x),u_3(x))\in\Rbb^3$ at every point $x\in\Omega$, 
then $v(x)$ is a matrix \begin{equation}
v(x)=\del{\dif u_1(x),\dif u_2(x), \dif u_3(x)}=\begin{pmatrix}
v_{11}(x) & v_{21}(x) & v_{31}(x) \\
v_{12}(x) & v_{22}(x) & v_{32}(x)
\end{pmatrix}
\end{equation}
where $v_{lj}(x):=\partial u_j(x)/\partial x_l$ at every $x=(x_1,x_2)$ for $l=1,2$ and $j=1,\dots.m$ ($m=3$ here).
As a result, the VTV of $u$ in \eqref{eq:VTV} simplifies to $\int_\Omega \|v(x)\|_\star \dif x$. 
If $u$ is not differentiable, then $v$ may become the weak gradients of $u$
(when $u\in W^{1,1}(\Omega;\Rbb^m)$) or even a Radon-Nikodym measure (when $u\in BV(\Omega;\Rbb^m)$).
However, in common practice using finite differences for numerical implementation
to approximate partial derivatives of functions,
such subtlety does not make much differences.
Therefore, we treat $v$ as the Jacobian of $u$ throughout the
rest of the paper and in numerical experiments.

Assume that we can derive the relation of the Jacobian $v$ and the original
observed data $b$ and form a data fidelity $H(v)$ of $v$, as an analogue of
data fidelity $h(u)$ of image $u$ (justification of this assumption 
will be presented in Section \ref{subsec:data_reform}). 
Then we can reformulate the VTV regularized 
image reconstruction problem \eqref{eq:VTVbased} about $u$ to a matrix $\star$-norm 
regularized inverse problem about $v$ as follows:
\begin{equation}\label{eq:L1based}
\min_v \alpha \|v\|_\star + H(v).
\end{equation}
Now we seek for reconstructing $v$ instead of $u$. 
This reformulation \eqref{eq:L1based} 
has two significant advantages compared to
the original formulation \eqref{eq:VTVbased}:
\begin{itemize}
\item It reduces to an $l_1$ type minimization \eqref{eq:L1based} and can be solved
effectively by accelerated gradient method.
For example, using Algorithm \ref{alg:edgeFISTA}
based on FISTA \cite{Beck:2009b}, we can attain
an optimal convergence rate of $O(1/k^2)$ to solve \eqref{eq:L1based},
where $k$ is the iteration number.
This is in sharp contrast to the best known $O(1/k)$ rate of primal-dual based 
methods (including the recent, successful PDHG and ADMM) 
for solving \eqref{eq:VTVbased}.
\item The per-iteration complexity of the $l_1$ type minimization \eqref{eq:L1based},
e.g., Algorithm \ref{alg:edgeFISTA}, is the same or even lower compared to
VTV regularized minimization \eqref{eq:VTVbased}. In particular,
closed-form solution of gradients in \eqref{eq:vsubp}, a matrix $\star$-norm
variant of soft-shrinkage, is widely available and cheap to compute.
\end{itemize}

\begin{algorithm}
\caption{Joint edge reconstruction of \eqref{eq:L1based} using FISTA}
\label{alg:edgeFISTA}
\begin{algorithmic}
\State \textbf{Input:} Initial $u_0$ and its Jacobian $v_0=D u_0$. Step size $\tau\leq 1/\|\nabla H\|$.
\State Set $k=0$, $t_0=1$, $w_0=v_0$, and iterate \eqref{eq:vsubp}--\eqref{eq:wsubp} below
until stopping criterion is met
\begin{align}
v_{k+1} & = \argmin_{v} \del[3]{\alpha \|v\|_\star + \frac{1}{2}\|v-w_k+\tau \nabla H(w_k)\|_F^2} \label{eq:vsubp}\\
t_{k+1} & = \frac{1+\sqrt{1+4t_k^2}}{2} \label{eq:tsubp}\\
w_{k+1} & = v_{k+1} + \frac{t_k-1}{t_{k+1}} (v_{k+1}-v_k) \label{eq:wsubp}
\end{align}  
\State \textbf{Output:} Reconstructed Jacobian (edge) $v \leftarrow v_{k+1}$.
\end{algorithmic}
\end{algorithm}

In the remainder of this section, we will answer the following three questions 
that are critical in the proposed joint image reconstruction method 
based on \eqref{eq:L1based}:
\begin{itemize}
\item In what situations/applications the fidelity $h(u)$ of image $u$ can be converted to fidelity $H(v)$ of Jacobian $v$? 
(Section \ref{subsec:data_reform})
\item How to obtain closed form solution of $v$-subproblem \eqref{eq:vsubp} for matrix $\star$-norms, particularly for the
three matrix norms in Section \ref{subsec:VTV}? (Section \ref{subsec:shrinkage})
\item How to reconstruct image $u$ using Jacobian $v$ obtained from Algorithm \ref{alg:edgeFISTA}? (Section \ref{subsec:recon_u})
\end{itemize}

\subsection{Formulating fidelity of Jacobian $v$}\label{subsec:data_reform}
In an extensive variety of imaging technologies, especially medical imaging, the image data 
are acquired in transform domains. Among those common transforms,
Fourier transform and Radon transform are the two most widely used ones in medical imaging.
In what follows, we show that the relation between an image $u$ and its (undersampled) data $f$
can be easily converted to that between the Jacobian matrix $v$ and $f$ in the
Fourier domains.
This allows us to derive the data fidelity of Jacobian $v$ and reconstruct it properly
for MRI reconstruction problems.
Similar idea can be applied to the case with Radon transforms,
subject to modification of data fidelity and noise distribution in a straightforward manner.
In this paper, we focus on the Fourier case only as our data are multi-contrast MRIs.

\textbf{Data transformation.}
Imaging technologies, such as magnetic resonance imaging (MRI) and radar imaging,
are based on Fourier transform of images. The image data are the Fourier coefficients of
image acquired in the Fourier domain (also known as the $k$-space in MRI community). 
The inverse problem of image reconstruction in such technologies often refers to
recovering image from partial (i.e., undersampled) Fourier data.

In discrete setting, let $u\in\Rbb^n$ denote the image to be reconstructed, 
$F\in\Rbb^{n\times n}$ the (discrete) Fourier transform matrix (hence a unitary matrix),
and $P\in\Rbb^{n\times n}$ the diagonal matrix with binary values ($0$ or $1$) as diagonal entries
to represent the undersampling pattern (also called mask in $k$-space),
and $f\in\Rbb^n$ be the observed partial data with $0$ at unsampled locations.
Therefore, the relation between the underlying image $u$ and observed partial data $f$
is given by $PFu=f$.

The gradient (partial derivatives) of an image can be regarded as a convolution.
The Fourier transform, on the other hand, is well known to convert a convolution to
simple point-wise multiplication in the transform domain.
In discrete settings, this simply means that $\hat{D}_i:=F D_i F^T$ is a diagonal
matrix where $D_i\in\mathbb{R}^{n\times n}$ is the discrete partial derivative (e.g., forward finite difference)
operator along the $x_i$ direction ($i=1,2$). This amounts to a straightforward formulation
for the data fidelity of $v$: 
Let $v_i$ be the partial derivative of $u$ in $x_i$ direction, i.e., $v_i=D_i u$ then we have
\begin{equation}\label{eq:data_trans}
PFv_i=PFD_iu=PFD_i F^T F u = P\hat{D}_iFu=\hat{D}_iPFu=\hat{D}_if
\end{equation}
where we used the facts that both $P$ and $\hat{D}_i$ are diagonal matrices 
so they can commute. Therefore, the data fidelity term $H(v)$ in \eqref{eq:L1based}
can be, for example, formulated as
\begin{equation}
H(v) = \frac{1}{2} \del[2]{ \| PFv_1 - \hat{D}_1 f\|^2 + \| PFv_2 - \hat{D}_2 f\|^2 }
\end{equation}
As long as the Fourier data is concerned, the data fidelity of the
corresponding gradient can be formulated in a similar way as above.

\textbf{Noise transformation.} We now consider the noise distribution in image reconstruction from Fourier data.
Suppose the noise is due to acquisition in Fourier transform domain
such that
\begin{equation}\label{eq:data_trans_noise}
f(\omega) = P \Fcal[u](\omega) + e(\omega)
\end{equation}
for each frequency value $\omega$ in Fourier domain. 
Here $\Fcal[u]$ is the Fourier transform of image $u$.
By multiplying $\Dhat_i(\omega)$ on both sides of \eqref{eq:data_trans_noise}, 
we see that this fidelity of $v_i$ 
can be obtained by
\[
\Dhat_i(\omega)f(\omega) = P \Fcal[v_i] (\omega) + \Dhat_i(\omega)e(\omega)\ .
\]
Suppose that $\Dhat_i(\omega)=A_i(\omega)+ \irm B_i(\omega)$ where $A_i(\omega)\in\Rbb$ and $B_i(\omega)\in\Rbb$ are real and imaginary parts
of $\Dhat_i(\omega)$ respectively, and $e(\omega)=a(\omega)+ \irm b(\omega)$ where $a(\omega)\sim N(0,\sigma_a^2)$
and $b(\omega)\sim N(0,\sigma_b^2)$ are independent for all $\omega$.
Then we know that the transformed noise $\Dhat_i(\omega)e(\omega)$ are independent for different $\omega$
and distributed as bivariate Gaussian 
\begin{equation}
\begin{pmatrix}
\text{real}(\Dhat_i(\omega)e(\omega))\\
\text{imag}(\Dhat_i(\omega)e(\omega))
\end{pmatrix}\sim
N\del{0,C_i(\omega) \Sigma C_i^T(\omega)}
\end{equation}
where $C_i(\omega)=[A_i(\omega),-B_i(\omega); B_i(\omega), A_i(\omega)]\in\Rbb^{2\times2}$ and $\Sigma = \text{diag}(\sigma_a^2,\sigma_b^2)\in\Rbb^{2\times2}$.
Therefore, we can readily obtain the maximum likelihood of $\Dhat_i(\omega)e(\omega)$ and hence
the fidelity of $v_i$.
In particular, if $\sigma_a=\sigma_b=\sigma$, then
$\text{real}(\Dhat_i(\omega)e(\omega))$ and $\text{imag}(\Dhat_i(\omega)e(\omega))$ are two i.i.d. $N(0,\sigma^2(A_i^2(\omega)+B_i^2(\omega)))$.
i.e., $N(0,\sigma^2|\Dhat_i(\omega)|^2)$.
In this case, we denote $\Psi_i=\text{diag}(1/|\Dhat_i(\omega)|^2)$, then 
data fidelity (i.e., negative log-likelihood)
of $v_i$ simply becomes $(1/2)\|P\Fcal v_i - \Dhat_i f\|_{\Psi_i}^2$,
as in contrast to the fidelity term $(1/2)\|P\Fcal u - f\|^2$ of $u$.
In the remainder of this paper, we assume
that the standard deviation of real and imaginary parts are both $\sigma$.
Furthermore, in numerical experiments, we can pre-compute $|\Dhat_i(\omega)|^2$
and only perform an additional point-wise division of $|\Dhat_i(\omega)|^2$
after computing $P\Fcal v_i-\Dhat_i f$ in each iteration.

\subsection{Closed form solution of matrix-valued shrinkage}\label{subsec:shrinkage}
As we can see, the algorithm step \eqref{eq:vsubp} calls for solution of
type
\begin{equation}\label{eq:mtxshrink}
\min_{X\in \Rbb^{2\times m}} \alpha \|X\|_\star + \frac{1}{2}\|X-B\|_F^2
\end{equation}
for specific matrix norm $\|\cdot\|_\star$ and given matrix $B\in\Rbb^{2\times m}$.
In what follows, we provide the close-form solutions of \eqref{eq:mtxshrink}
when the matrix $\star$-norm is Frobenius, induced 2-norm or nuclear norm, as mentioned in Section \ref{subsec:VTV}.
The derivations are provided in Appendix \ref{apd:shrinkage}.
\begin{itemize}
\item \textbf{Frobenius norm.}
This is the simplest case since $\|X\|_F$ treats $X$ as a vector in $\Rbb^{2m}$ in \eqref{eq:mtxshrink},
for which the shrinkage has close-form solution. More specifically, the solution of \eqref{eq:mtxshrink} is
\begin{equation}\label{eq:Fro_sol}
X^*=\max(\|B\|_F-\alpha,0)\frac{B}{\|B\|_F}\quad .
\end{equation}

\item \textbf{Induced 2-norm.}
This is advocated the vectorial TV \cite{Goldluecke:2010a}, but now we can provide a close-form solution
of \eqref{eq:mtxshrink} as
\begin{equation}\label{eq:2norm_sol}
X^*=B-\alpha \xi\eta^T
\end{equation}
where $\xi$ and $\eta$ are the left and right singular vectors corresponding to the largest singular value of $B$.

\item \textbf{Nuclear norm.}
This norm promotes low rank and yields a close-form solution of \eqref{eq:mtxshrink} as
\begin{equation}\label{eq:nuclear_sol}
X^*=U \max(\Sigma-\alpha,0) V^T
\end{equation}
where $(U,\Sigma,V)$ is the singular value decomposition (SVD) of $B$.
\end{itemize}
The computation of \eqref{eq:Fro_sol} is essential shrinkage of vector and hence very cheap.
The computations of \eqref{eq:2norm_sol} and \eqref{eq:nuclear_sol} involve
(reduced) SVD, however, explicit formula also exists as the matrices have tiny size of $2$-by-$m$
($3$-by-$m$ if images are 3D), where $m$ is the number of image channels/modalities.
\smallskip

\noindent \textbf{Remark.}
It is worth noting that the computation of \eqref{eq:mtxshrink} is carried out at every pixel
independently of others in each iterations. This allows straightforward parallel computing 
which can further reduce real-world computation time.


\subsection{Reconstruct image from gradients}\label{subsec:recon_u}
Once the Jacobian $v$ is reconstructed from data, the final step is to
resemble the image $u$ from $v$. Since this step is performed for each
modality, the problem reduces to reconstruction of a scalar-valued image $u$ from its gradient $v=(v_1,v_2)$.
In \cite{Patel:2012a,Sakhaee:2015a}, the image $u$ is 
reconstructed by solving the Poisson equation $\Delta u=\mathrm{div}(v)=-(D_1^Tv_1+D_2^Tv_2)$
since $v_1=D_1u$ and $v_2=D_2u$ are the partial derivatives of $u$.
The boundary condition of this Poisson equation can be either Dirichlet or Neumann
depending on the property of imaging modality. 
In medical imaging applications, such as MRI, CT, and PET, the boundary condition
is simply $0$ since it often is just background near image boundary $\partial\Omega$.

Numerically, it is more straightforward to recover $u$ by solving 
the following minimization
\begin{equation}\label{eq:u_recon}
\min_{u} \|D_1u-v_1\|^2 + \|D_2u-v_2\|^2 + \beta h(u)
\end{equation}
with some parameter $\beta>0$ to weight the data fidelity $h(u)$ of $u$.
The solution is easy to compute since the objective is smooth, and often times
closed-form solution may exist. For example, in the MRI case 
where $h(u)=\frac{1}{2}\|PFu-f\|^2$, the solution of \eqref{eq:u_recon} is
given by
\begin{equation}\label{eq:usol}
u = \Fcal^T\sbr[3]{ \del[2]{\Dhat_1^T\Dhat_1+\Dhat_2^T\Dhat_2+\beta P^TP}^{-1} 
\del[2]{\Dhat_1^T\Fcal v_1+\Dhat_2^T\Fcal v_2+\beta P^Tf} }
\end{equation}
where $\Dhat_1^T\Dhat_1+\Dhat_2^T\Dhat_2+\beta P^TP$ is diagonal and hence
trivial to invert. The main computations are just few Fourier transforms.

In fact, it seems often sufficient to retrieve the base intensity of $u$ 
from $h(u)$ to reconstruct $u$ from $v$, and hence the result is not sensitive to $\beta$ when solving \eqref{eq:u_recon}.
It is also worth pointing out that, by setting $\beta=0$, the minimization
\eqref{eq:u_recon} is just least squares whose normal equation is the Poisson equation mentioned above.

To summarize, we propose the two-step Algorithm \ref{alg:twostep}
which restores image edge $v$ and resembles image $u$ for multi-contrast MRI reconstruction.
The two steps are each executed only once (no iteration). Step 1 itself requires iterations
which converge very quickly with rate $O(1/k^2)$ where $k$ is iteration number.
Step 2 has closed form solution \eqref{eq:usol} for multi-contrast MRI reconstruction.
\begin{algorithm}
\caption{Two-step edge-based reconstruction for multi-contrast MRI}
\label{alg:twostep}
\begin{algorithmic}
\State \textbf{Input:} Multi-contrast partial Fourier data $f=(f_1,\dots,f_m)$, mask $P$. Initial $u_0$.
\State \textbf{Step 1:} Jointly reconstruct Jacobian $v=(v_1,\dots,v_m)$ using Algorithm \ref{alg:edgeFISTA}.
\State \textbf{Step 2:} Reconstruct image $u_j$ of contrast $j$ using $v_j$ and $f_j$ by \eqref{eq:usol} for $j=1,\dots,m$.
\State \textbf{Output:} Multi-contrast image $u=(u_1,\dots,u_m)$.
\end{algorithmic}
\end{algorithm}

\section{Numerical Results}
\label{sec:results}

In this section, we conduct a series of numerical experiments on synthetic and real 
multi-contrast MRI datasets using Algorithm \ref{alg:twostep} (for short, we call it
\textbf{ER}, standing for Edge-based Reconstruction, without adaptive weighting,
and \textbf{ER-weighted} for the one with adaptive weighting $\Psi$ in Section \ref{subsec:data_reform}).
For comparison, we also obtained implementation of two state-of-the-arts method for
joint image reconstruction: multi-contrast MRI method Bayesian CS 
(\textbf{BCS})\footnote{BCS code: \url{http://martinos.org/~berkin/Bayesian_CS_Prior.zip}} \cite{Bilgic:2011a}, 
and Fast Composite Splitting Algorithm for multi-contrast MRI (\textbf{FCSA-MT})\footnote{FCSA-MT code: \url{http://ranger.uta.edu/~huang/codes/Code_MRI_MT.zip}} 
\cite{Huang:2014a},
which are both publicly available online. 
We tune the parameters of each methods so they can perform nearly optimally.
In particular, FCSA-MT seems to be sensitive to noise level and we need to tune 
a very different parameter for different $\sigma$.
Specifically, we set the $l_1$ weight in BCS as $10^{-8}$ for both $\sigma=4$ and $10$ cases;
the TV term and wavelet term weights in FCSA-MT are set to $0.010$ and $0.035$ (resp.) for $\sigma=4$ case,
and $2.50$ and $2.50$ (resp.) for $\sigma=10$ case;
the $l_1$ weight of ER/ER-weighted is set to $10^{-3}$, and the fidelity weight for \eqref{eq:usol} to obtain
image from edges is set to $\beta=0.001$ for both $\sigma=4$ and $10$ cases.

The image datasets we used are a 2D multi-contrast Shepp-Logan phantom (size $256\times256$), 
a simulated multi-contrast brain image (size $256\times256$) obtained from 
BrainWeb\footnote{BrainWeb: \url{http://brainweb.bic.mni.mcgill.ca/brainweb/}}, and 
a 2D in-vivo brain image (size $217\times181$) included in the FCSA-MT code. 
All images contain three contrasts. In particular,
the three contrasts in BrainWeb and in-vivo brain datasets represent the T1, T2, and PD images.
Undersampling pattern
and ratio are presented below with results.
The experiments are performed in Matlab computing environment 
on a Mac OS with Intel i7 CPU and 16GB of memory. 
Gaussian white noise with standard deviation $\sigma=4$ and $10$
were added to the k-space for the simulated data. We 
full Fourier data to obtain reference (ground truth) image $u^*$. 
The relative (L2) error of reconstruction $u_j$ for modality $j$ is defined by
$\|u_j-u_j^*\|/\|u_j^*\|$ and used to measure error of $u_j$ quantitatively for $j=1,\dots,m$.

\subsection{Comparison of different matrix norms}
\begin{figure}[t!]
\centering
\begin{minipage}[b]{0.43\textwidth}
\includegraphics[width=\textwidth]{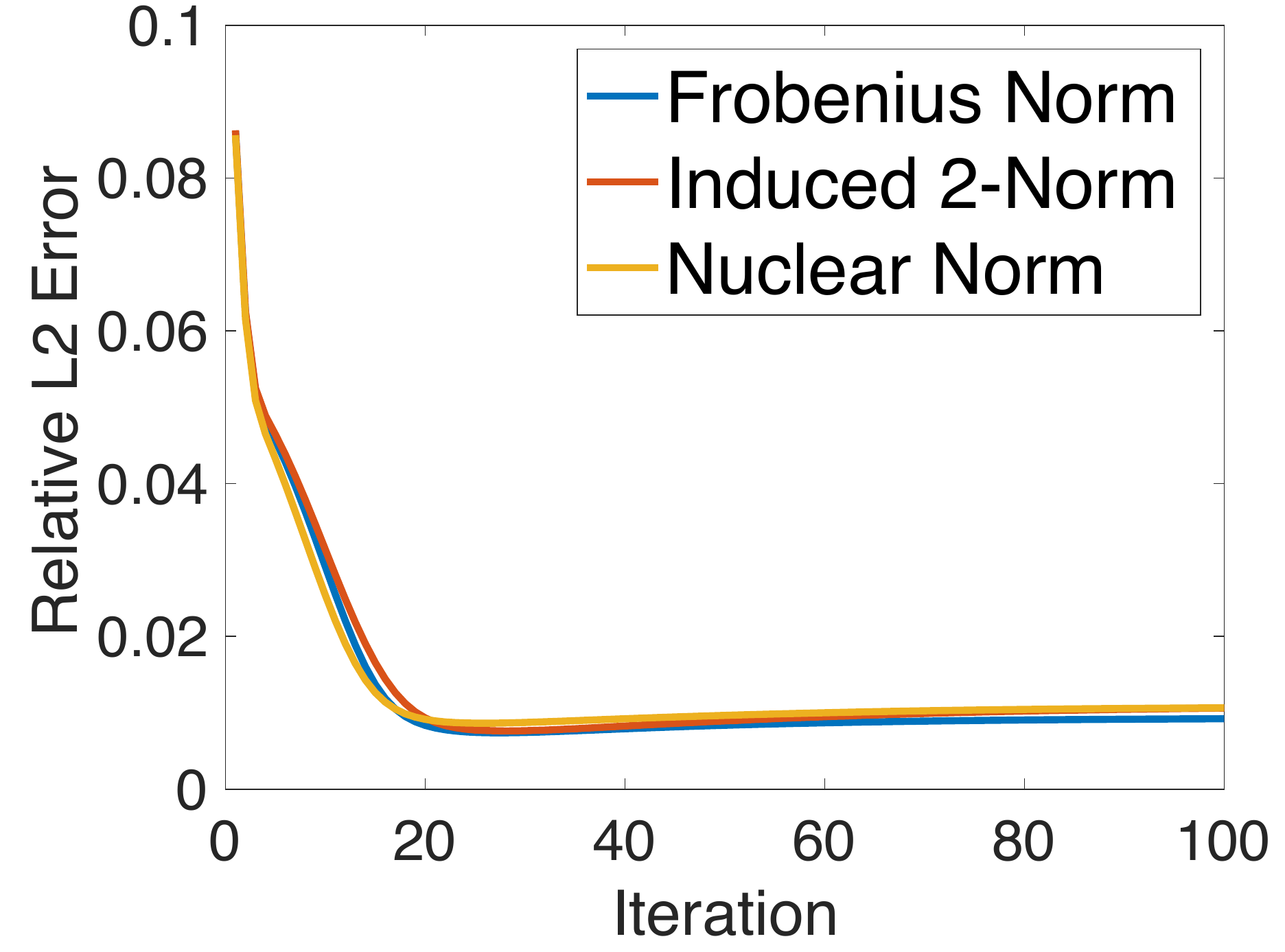}
\end{minipage}
\begin{minipage}[b]{0.41\textwidth}
\includegraphics[width=\textwidth]{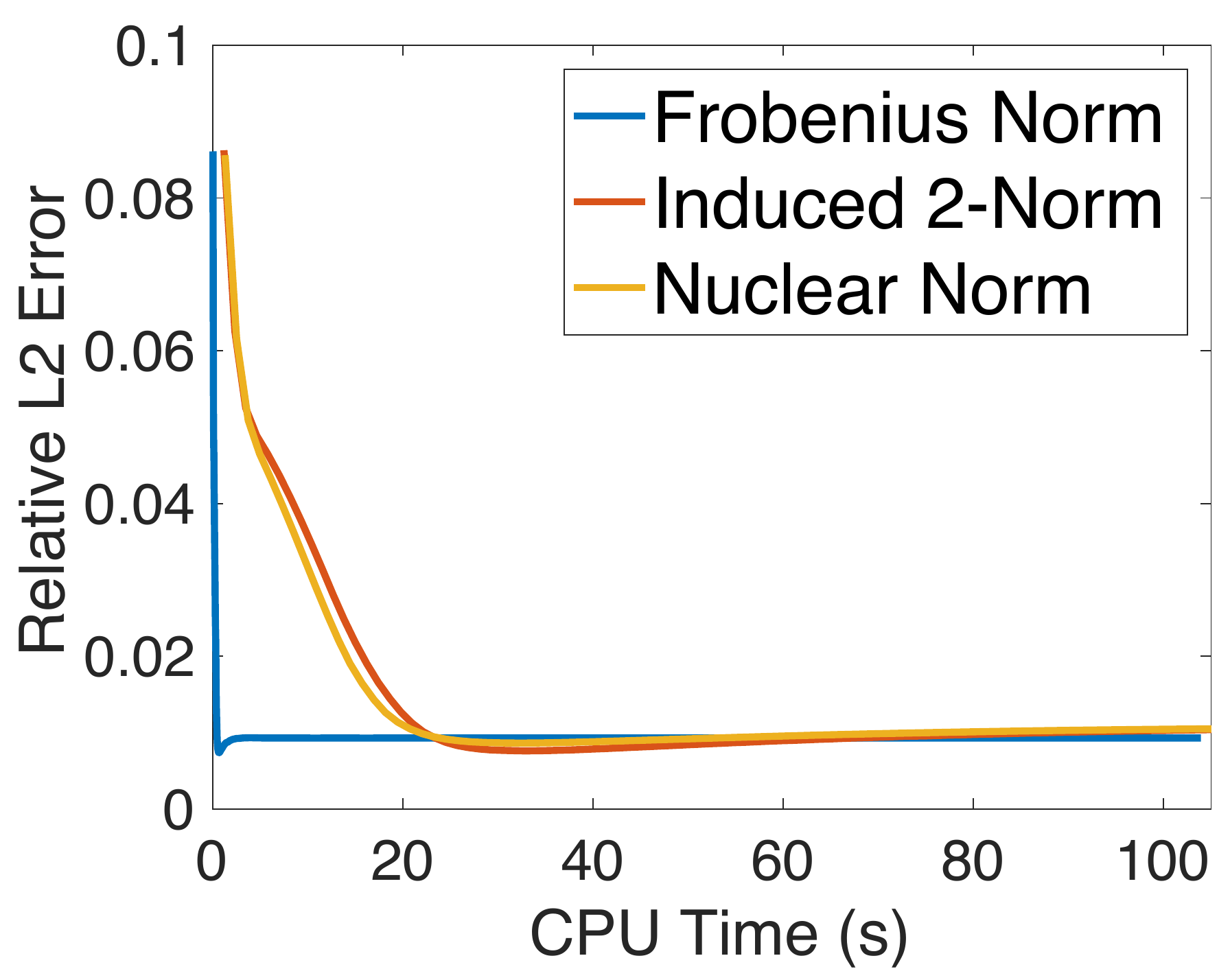}
\end{minipage}
\caption{Comparisons of Frobenius norm, induced 2-norm, and nuclear norm for Algorithm \ref{alg:twostep} 
(ER) on the BrainWeb image without noise. (a) Relative error vs number of iterations; (b) Relative error vs CPU time.}
\label{fig:norms}
\end{figure}

Our first test is on the performance of three matrix norms,
namely Frobenius, Induced 2-norm, and Nuclear norm, when used in Algorithm \ref{alg:twostep}.
These norms correspond to three commonly used vectorial TV regularization in multi-modal/channel/contrast
joint image reconstruction. All three norms yield closed-form solutions for subproblem \eqref{eq:vsubp}
as we showed in Section \ref{subsec:shrinkage}.
We apply the proposed Algorithm \ref{alg:twostep} with these three norms to
different image data and noise level combinations, and observe very similar
performance. For demonstration, we show a typical result using BrainWeb data with no noise 
in left panel of Figure \ref{fig:norms}.
As we can see, all norms yield very similar accuracy in terms of reconstruction error.

In terms of real-world computational cost, however, 
Frobenius norm treats matrix as vector and hence the shrinkage \eqref{eq:Fro_sol}
is very cheap to compute, while both induced 2-norm and nuclear norm 
require (reduced) SVD in \eqref{eq:2norm_sol} and \eqref{eq:nuclear_sol} 
that can be much slower
despite that the matrices all have small size ($2\times m$ for 2D images with $m$ modalities/contrasts).
We show the same trajectory of relative errors but versus CPU time
in the right panel of Figure \ref{fig:norms}, and it appears that Frobenius norm is the most 
cost-effective choice for this test. In the remainder of this section, we only use Frobenius norm
(corresponding to the specific vectorial TV norm of form $\int_\Omega \|Du(x)\|_F\dif x$)
for the proposed Algorithm \ref{alg:twostep},
but not the induced 2-norm (corresponding to $\int_\Omega \|Du(x)\|_2\dif x$) and
nuclear norm (corresponding to $\int_\Omega \|Du(x)\|_*\dif x$).

\subsection{Comparisons on multi-contrast MRI datasets}

\begin{figure}[t]
\centering
\includegraphics[width=.3\textwidth]{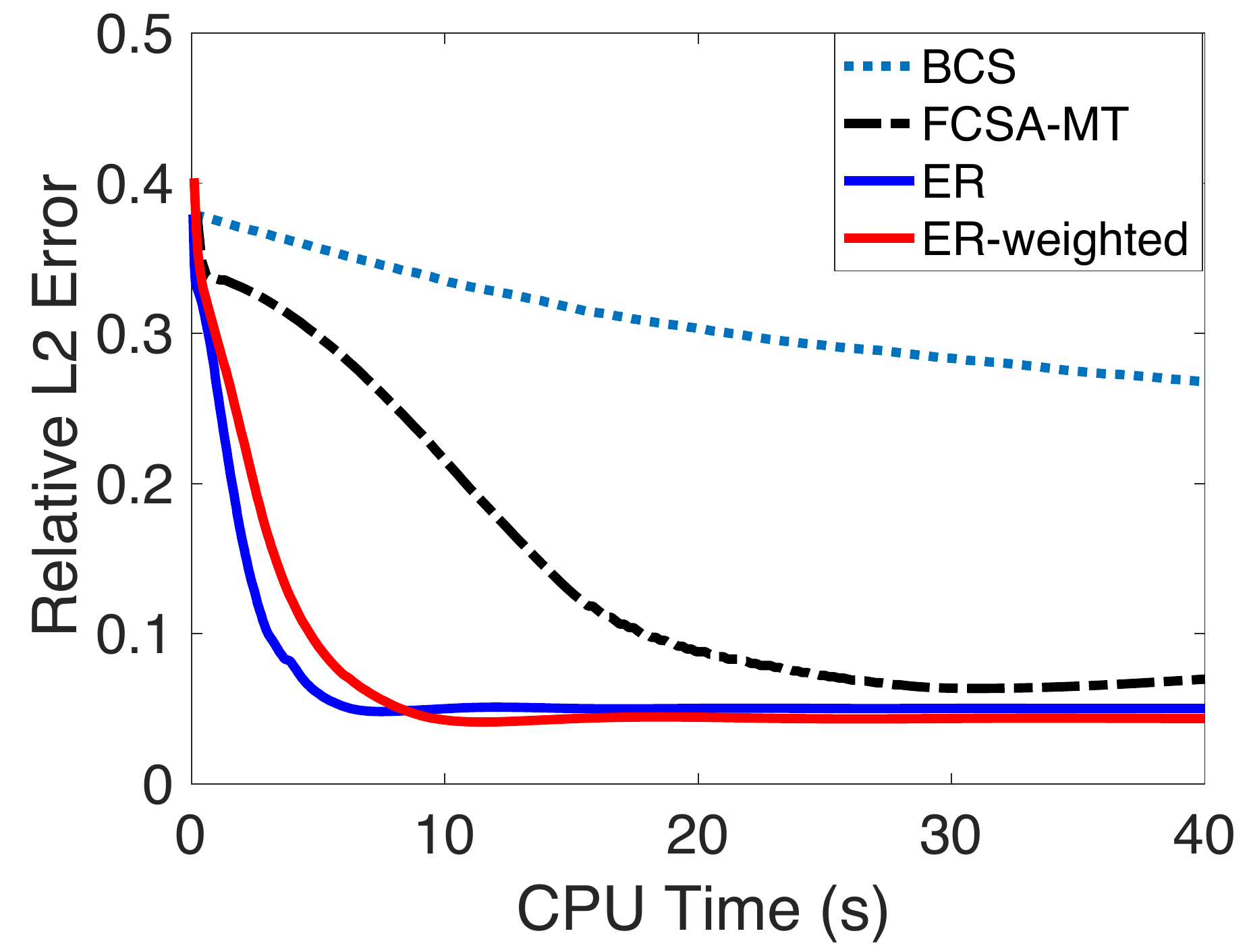}
\includegraphics[width=.3\textwidth]{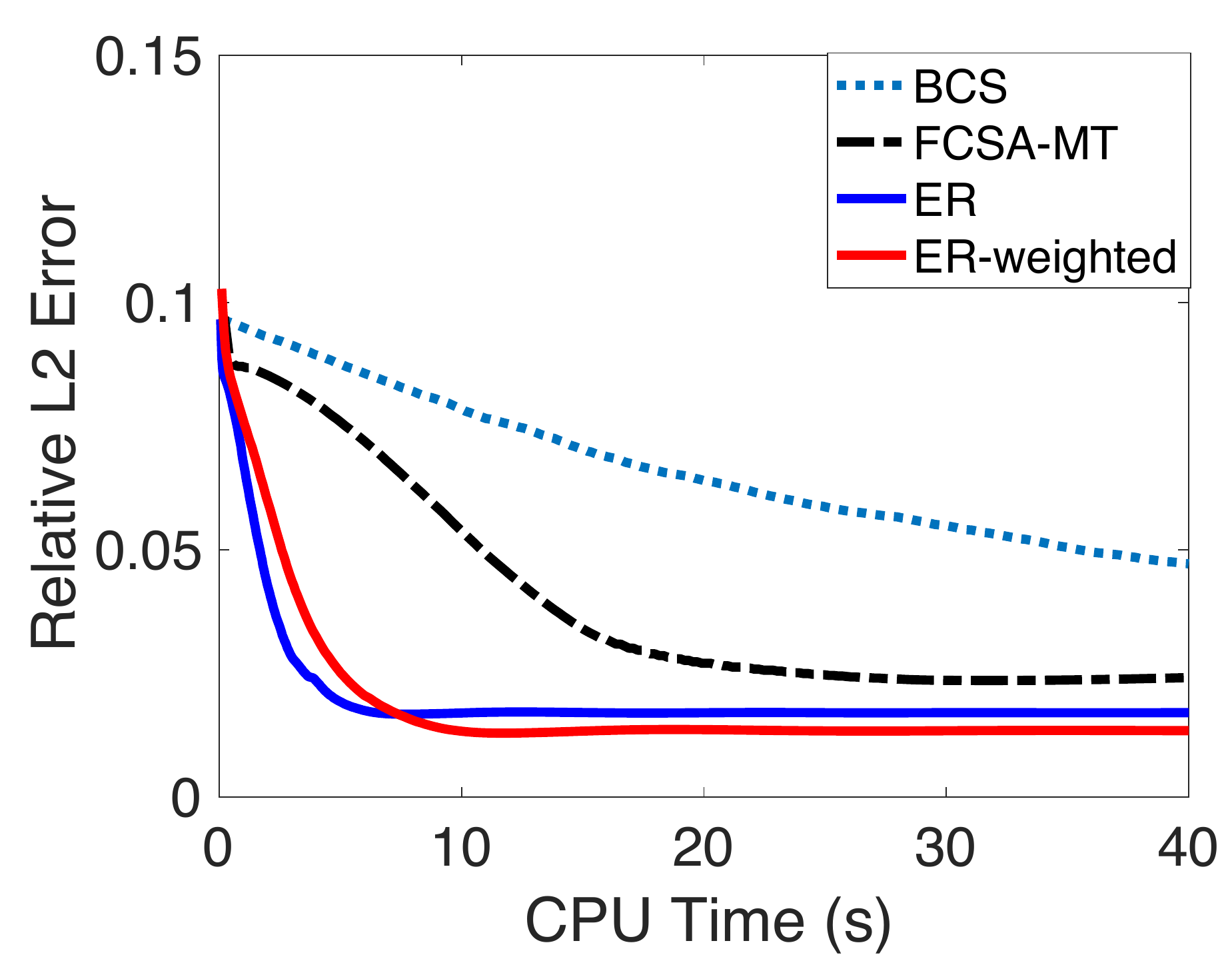}
\includegraphics[width=.3\textwidth]{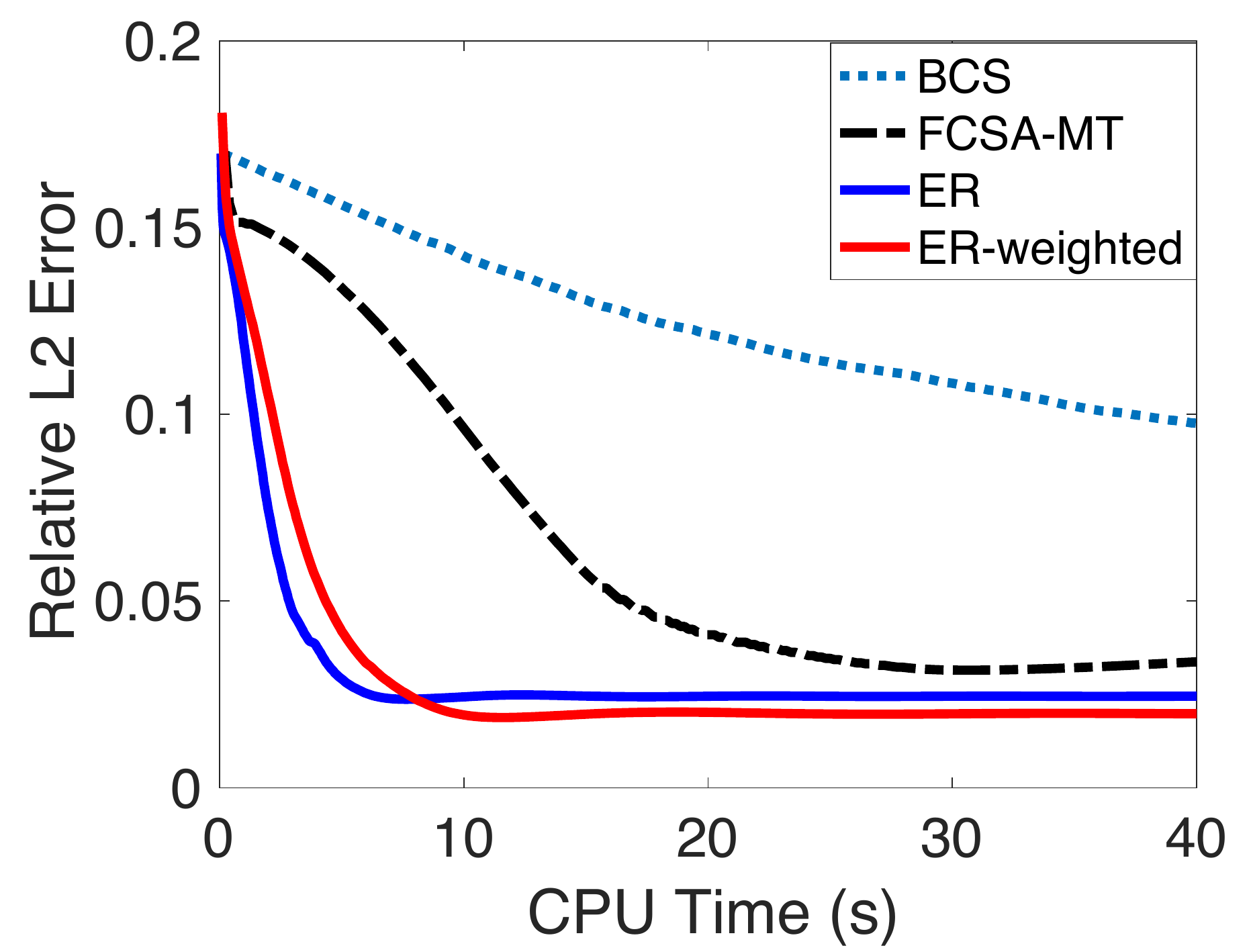}
\includegraphics[width=.3\textwidth]{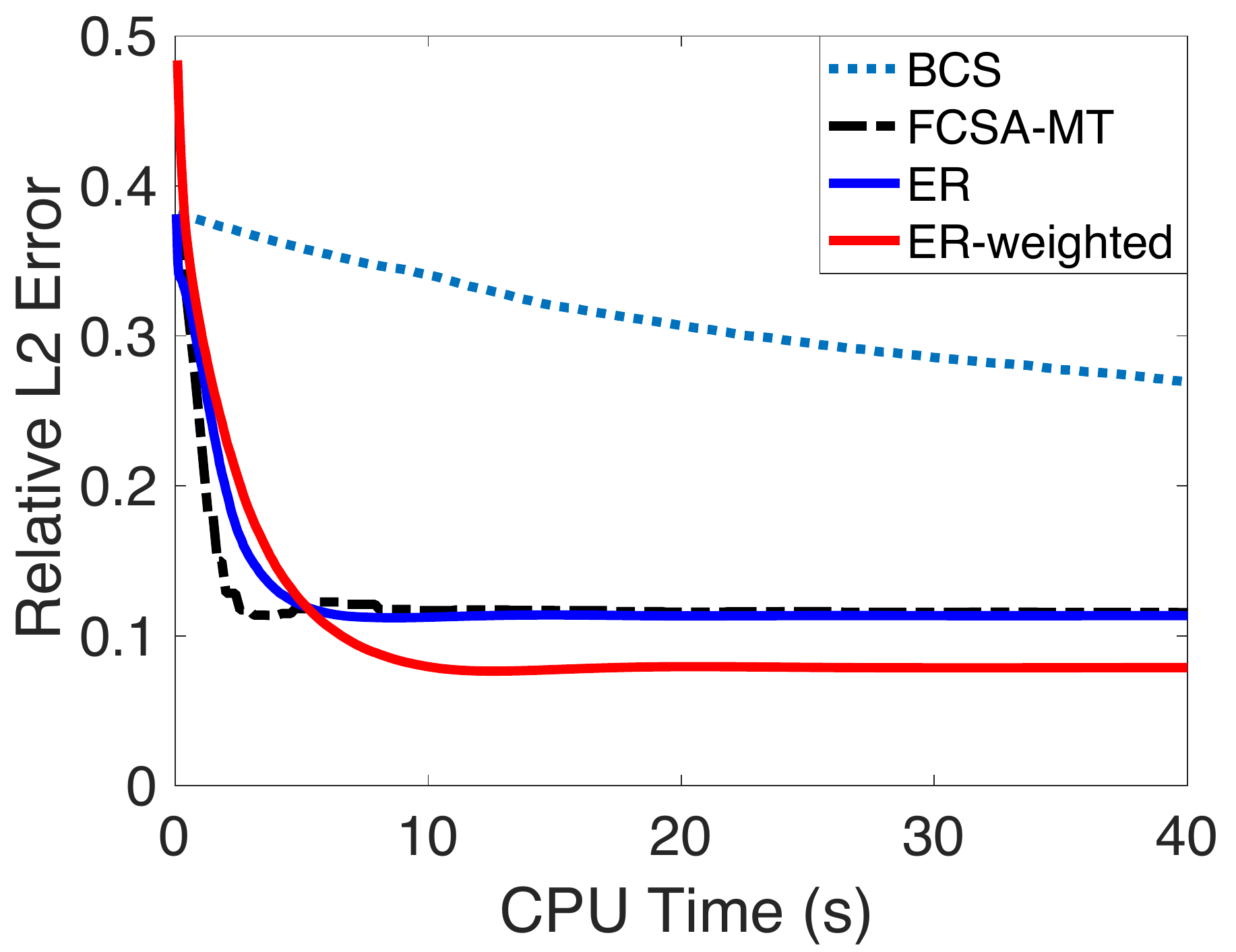}
\includegraphics[width=.3\textwidth]{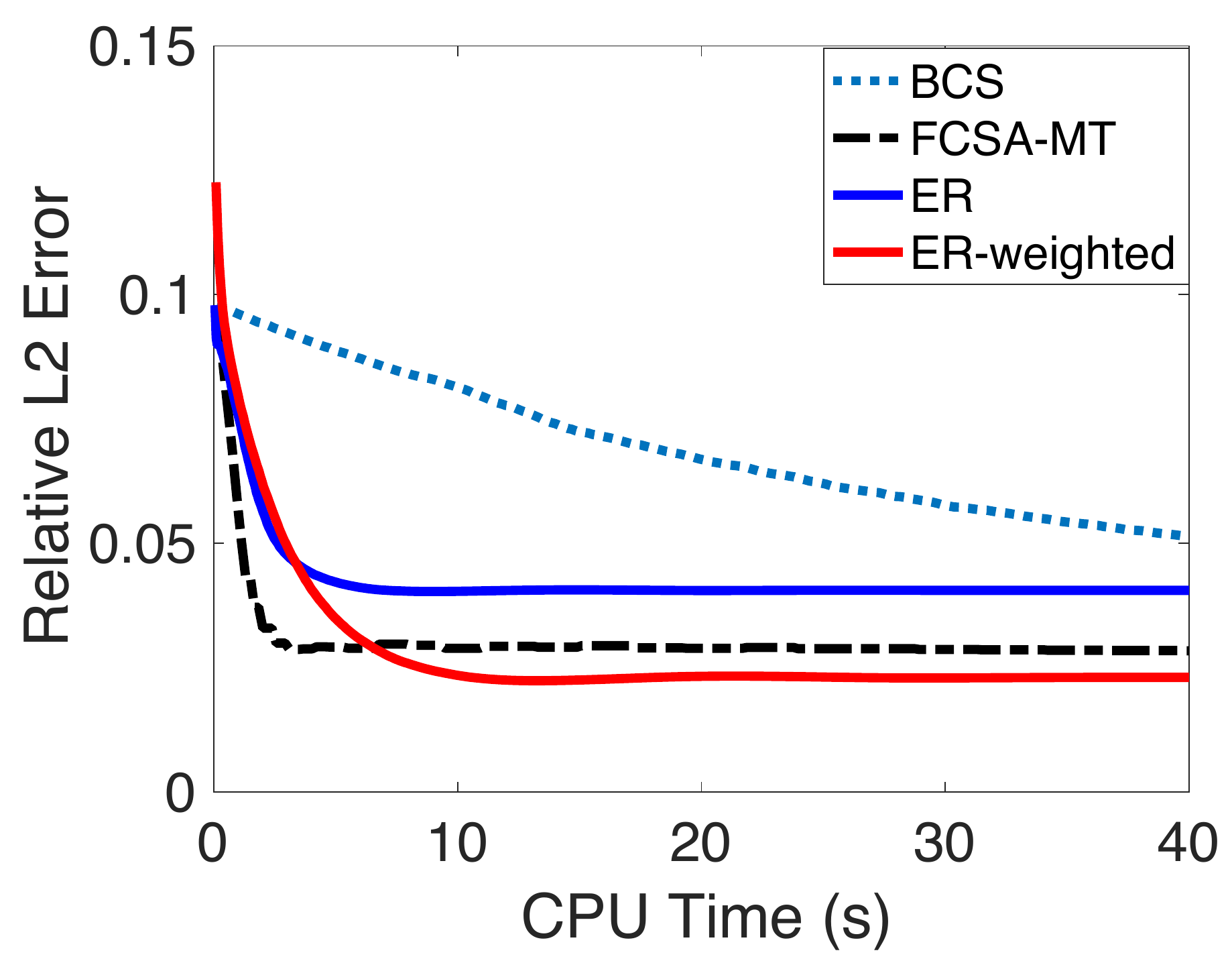}
\includegraphics[width=.3\textwidth]{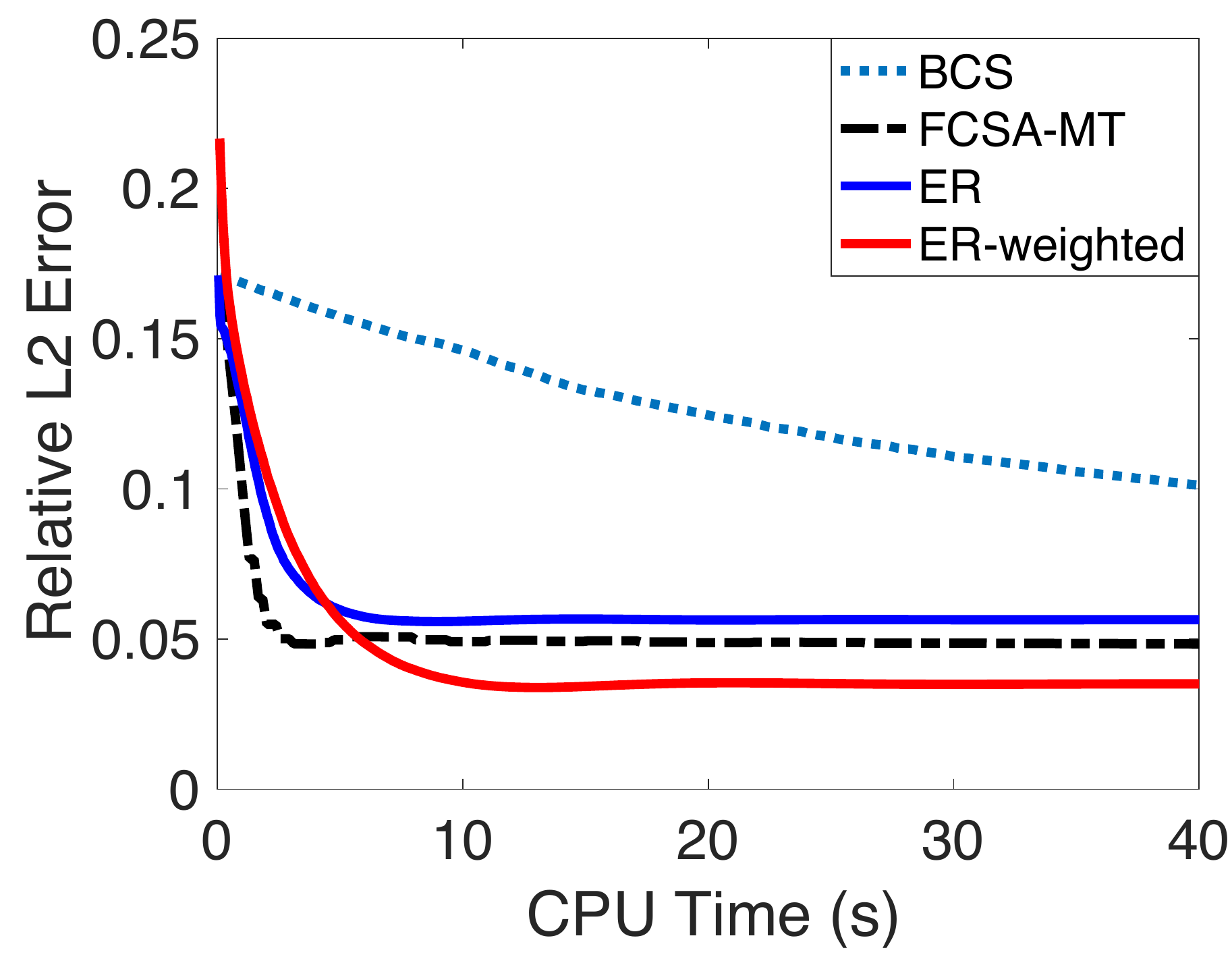}
\caption{Relative error vs CPU time of the comparison algorithms 
on Shepp-Logan phantom image with noise standard deviation $\sigma=4$ (top row)
and $\sigma=10$ (bottom row).} 
\label{fig:phantom_sigma_error_time}
\end{figure}

Now we conduct comparison of the proposed algorithm with BSC and FCSA-MT on multi-contrast MRI datasets.
We first test the comparison algorithms on a multi-contrast Shepp-Logan phantom image.
To demonstrate robustness of the proposed algorithm, 
we use a radial mask of undersampling ratio $13.5\%$ 
in Fourier space, but add white complex-valued Gaussian noise with standard deviation $\sigma=4,10$ (for both
real and imaginary parts of the noise) to the undersampled Fourier data.
Then we apply all comparison algorithm, namely BCS, FSCA-MT, and proposed ER and ER-weighted
to the undersampled noisy Fourier data to reconstruct the multi-contrast image.
For each method, we record the progress of relative error and show the relative error vs CPU time 
in Figure \ref{fig:phantom_sigma_error_time}. 
From these plots, we can see that the proposed Algorithm \ref{alg:twostep} (ER/ER-weighted) 
achieves higher efficiency
as they reduce reconstruction error faster than the comparison algorithms BCS and FCSA-MT.
Moreover, the proposed algorithm with the weight $\Psi$ (ER-weighted) can 
incorporate the transformed noise distribution and further improve reconstruction accuracy.
We observe that ER-weighted sometimes performs slightly slower than ER
in terms of CPU time (although always faster in terms of iteration number which is not shown here),
since the ER-weighted requires an additional point-wise division using $\Psi$.
We believe that this is an issue that can be resolved by further optimizing code, because complexity wise
the additional computation required by ER-weighted is rather low compared to other operations (such as
Fourier transforms) in ER/ER-weighted.

\begin{figure}[h!]
\centering
\includegraphics[width=.8\textwidth]{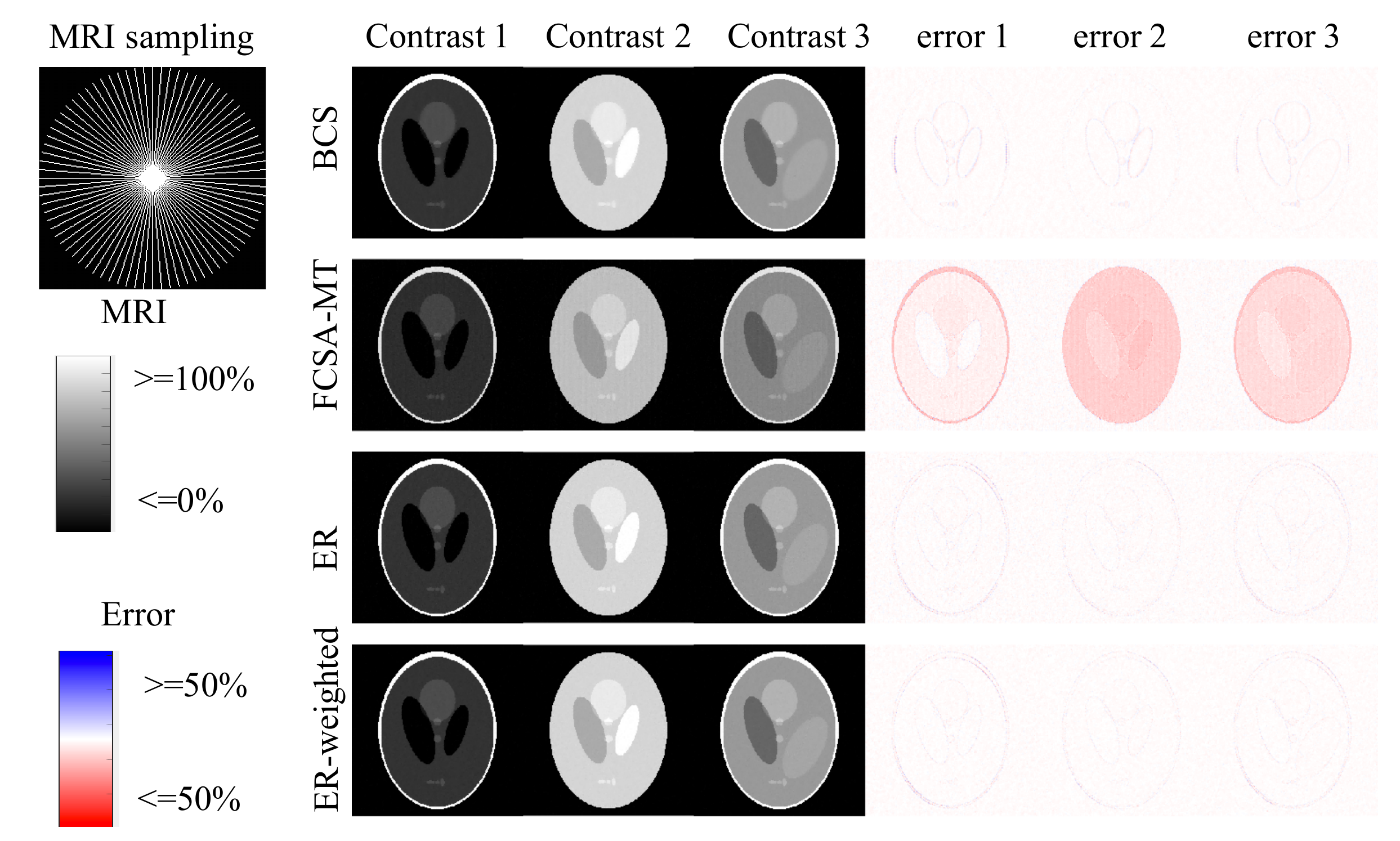}
\caption{Undersampling mask, reconstructed image, and error image for the Shepp-Logan image 
with noise level $\sigma=4$.} 
\label{fig:phantom_sigma4_visual}
\end{figure}

Besides the plot of relative error vs CPU time, we also show the final reconstruction images
using these methods in Figure \ref{fig:phantom_sigma4_visual} for the noise level $\sigma=4$ case.
In Figure \ref{fig:phantom_sigma4_visual}, the radial mask of sampling ratio (13.5\%)
is show on upper left corner (white pixels indicate sampled location).
The final reconstructed multi-contrast images $u=(u_1,u_2,u_3)$ by the four comparison algorithm
are shown in the middle column with algorithm name on left. On the right column,
we also show their corresponding error image $u_j-u_j^*$ where $u^*=(u_1^*,u_2^*,u_3^*)$ is the 
ground truth multi-contrast image obtained using full Fourier data.
As we can see, the reconstruction error of proposed algorithm is very small compared other others,
especially showing less obvious error on edges. This demonstrates the
improve reconstruction accuracy using our method. 


\begin{figure}[t]
\centering
\includegraphics[width=.32\textwidth]{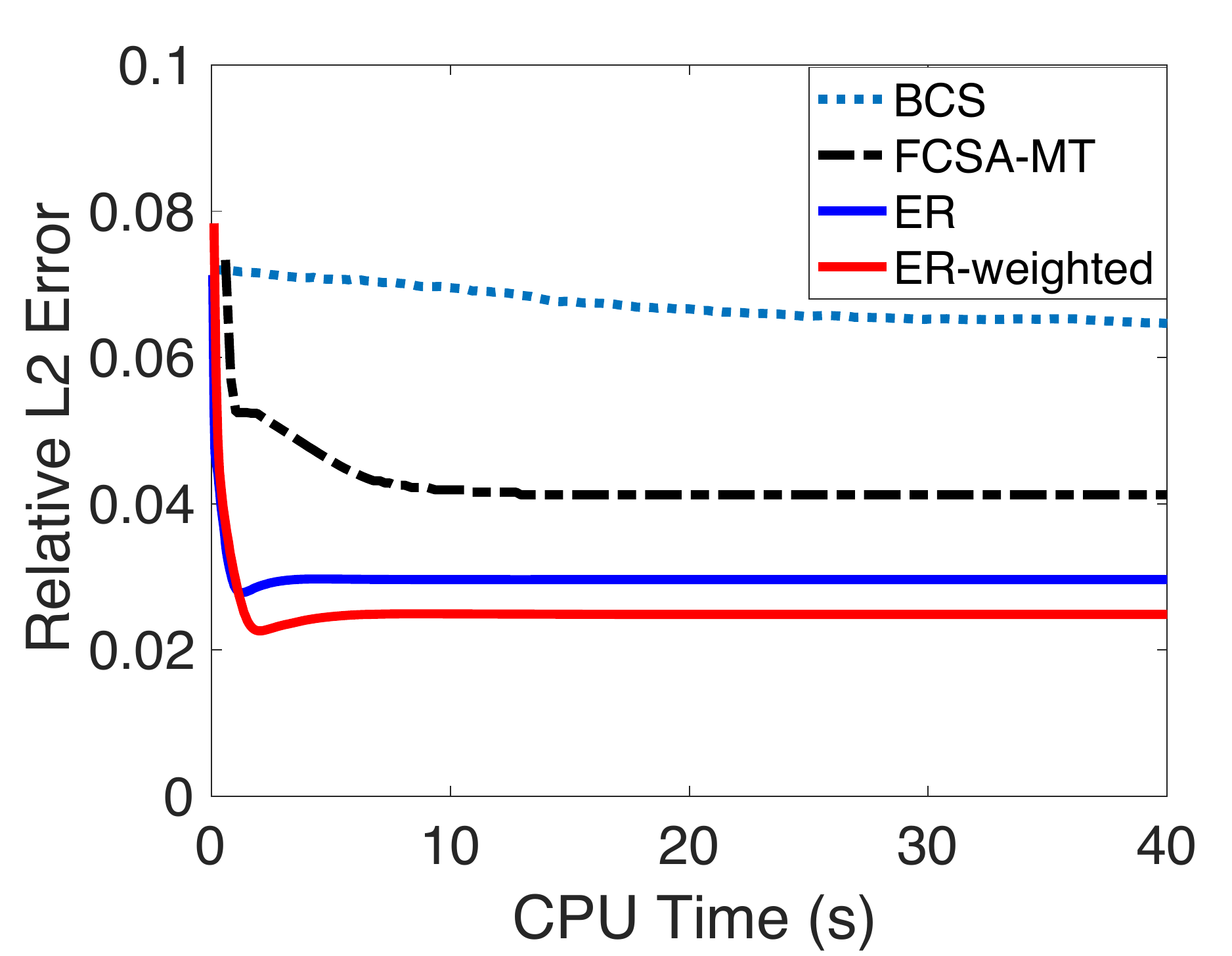}
\includegraphics[width=.32\textwidth]{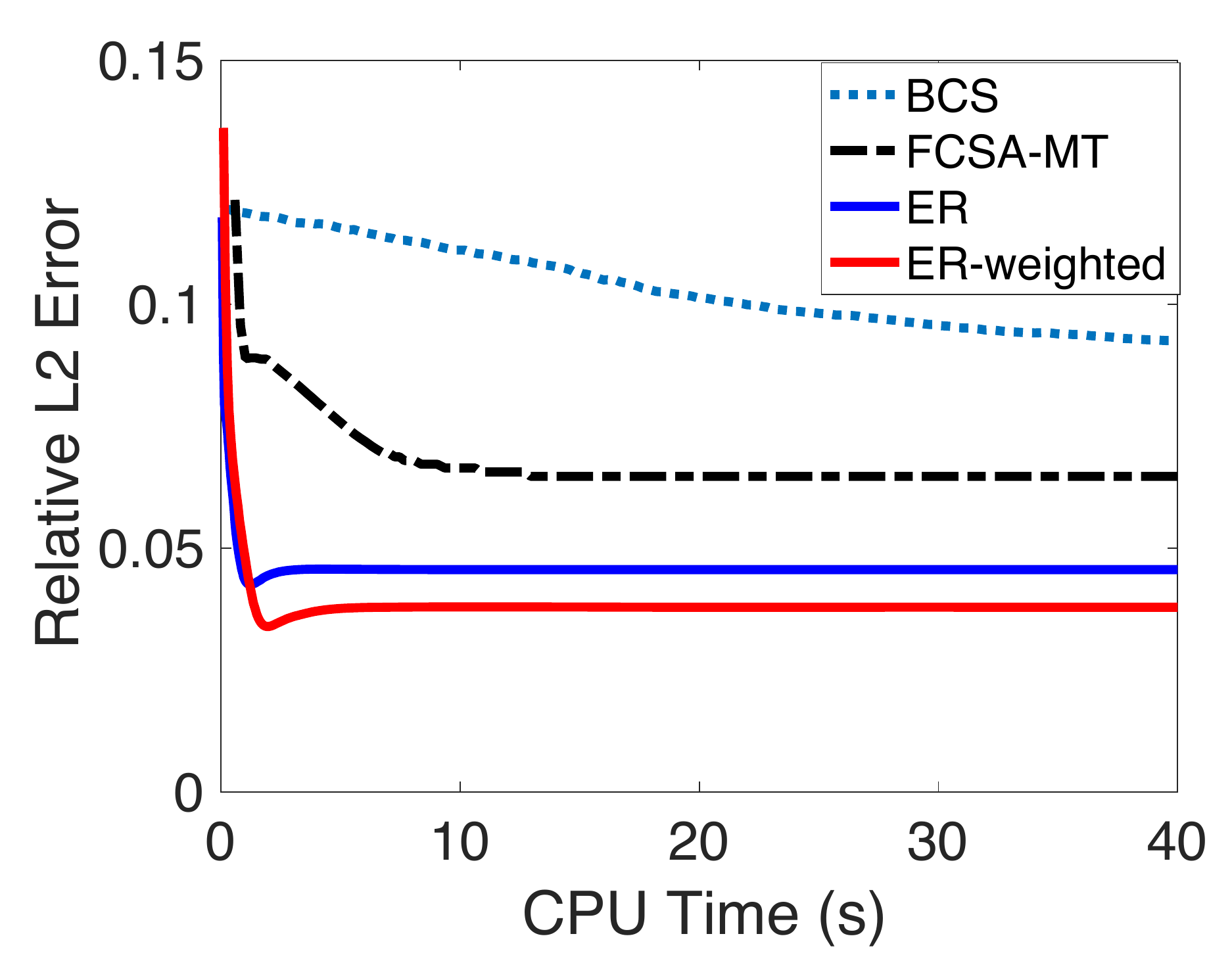}
\includegraphics[width=.32\textwidth]{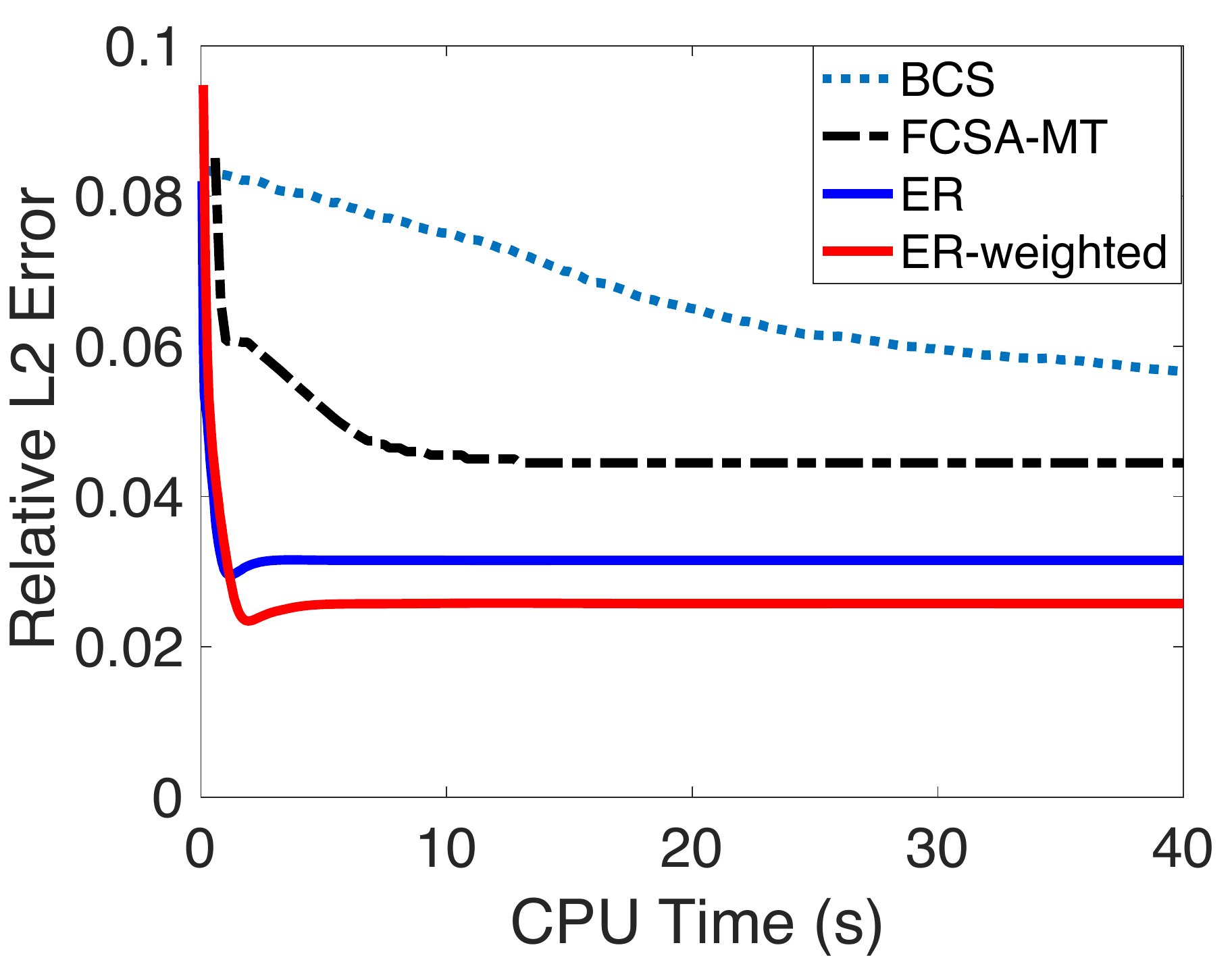}
\includegraphics[width=.32\textwidth]{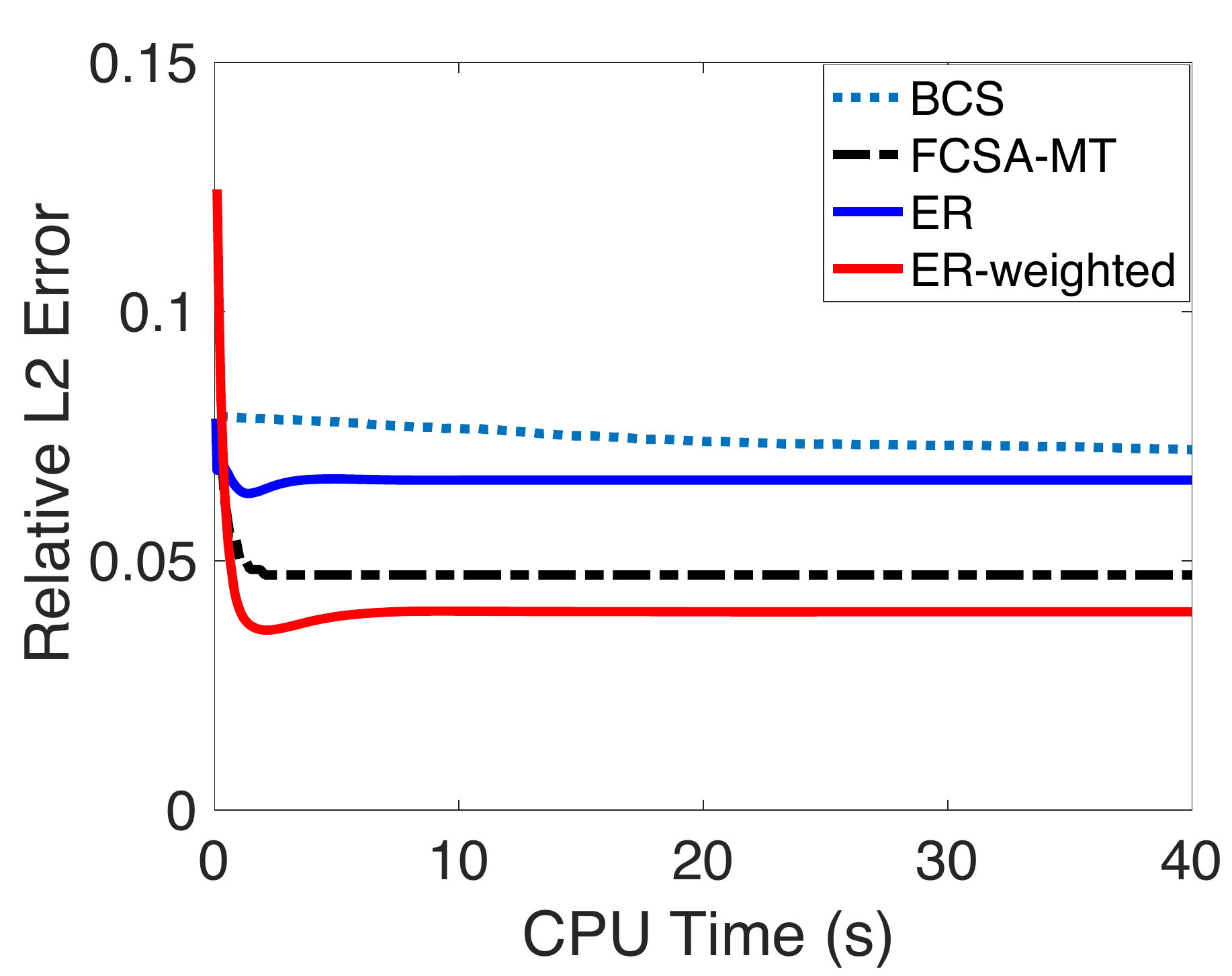}
\includegraphics[width=.32\textwidth]{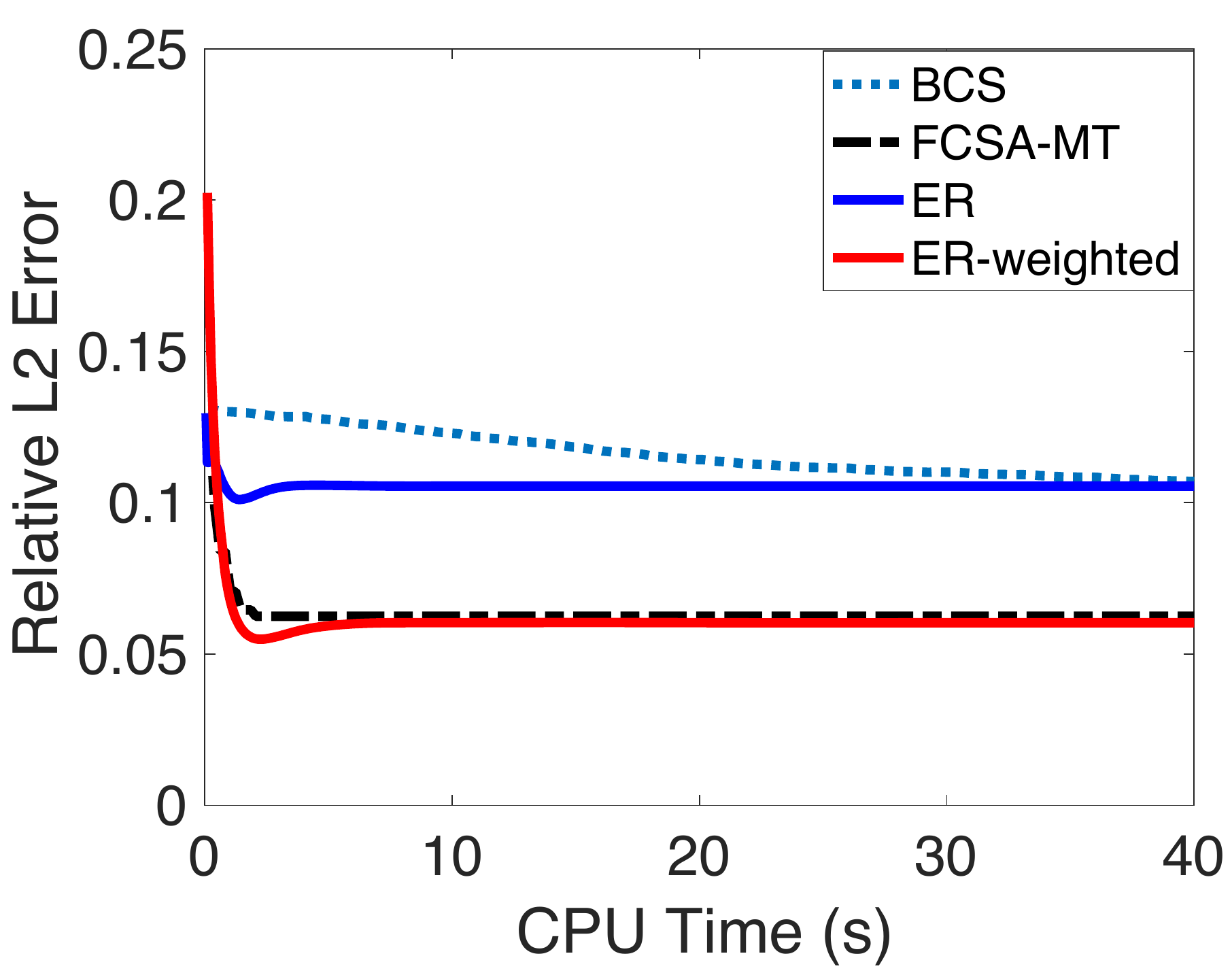}
\includegraphics[width=.32\textwidth]{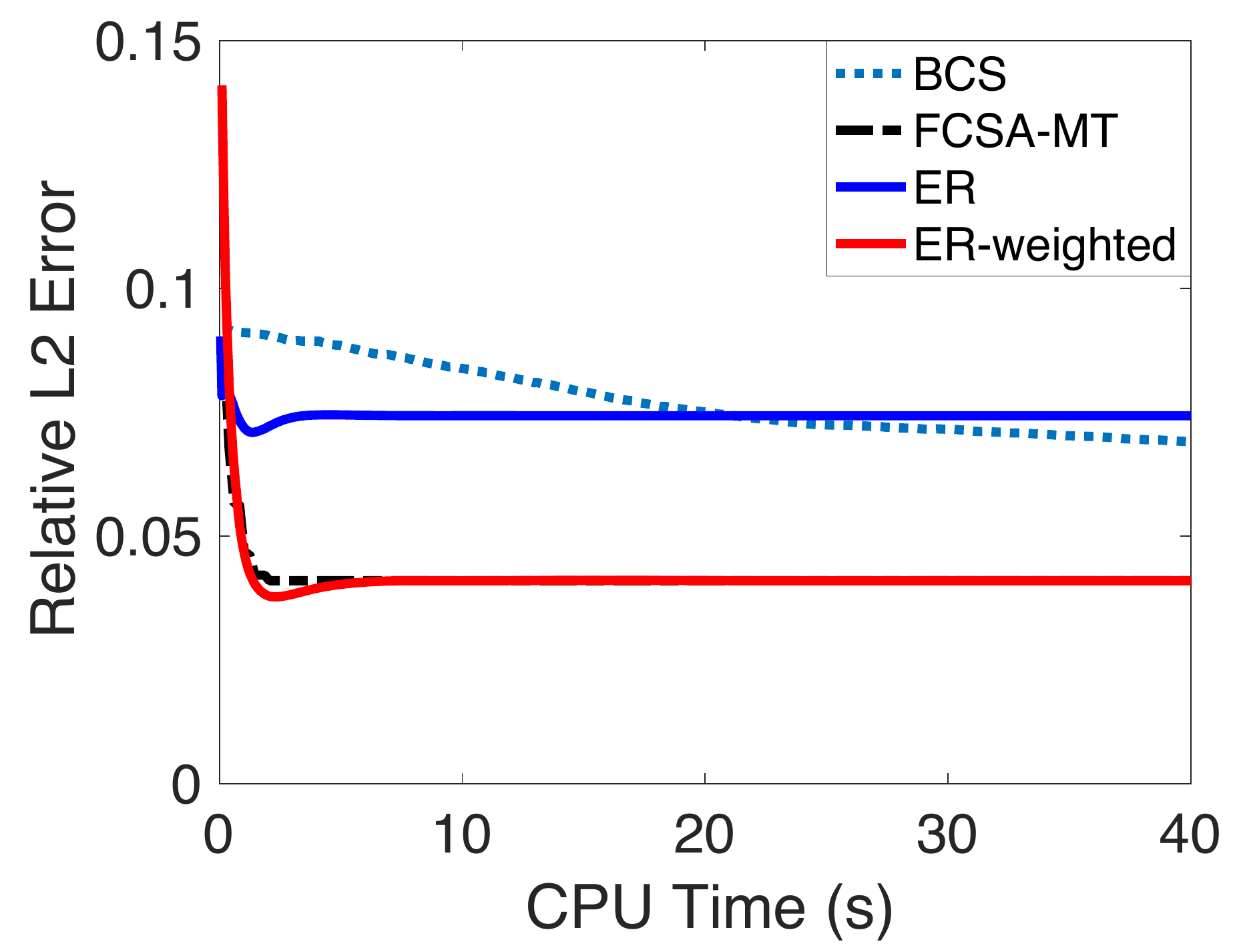}
\caption{Relative error vs CPU time of the comparison algorithms 
on BrainWeb image with noise standard deviation $\sigma=4$ (top row)
and $\sigma=10$ (bottom row).} 
\label{fig:brainweb_sigma_error_time}
\end{figure}

\begin{figure}[h!]
\centering
\includegraphics[width=.9\textwidth]{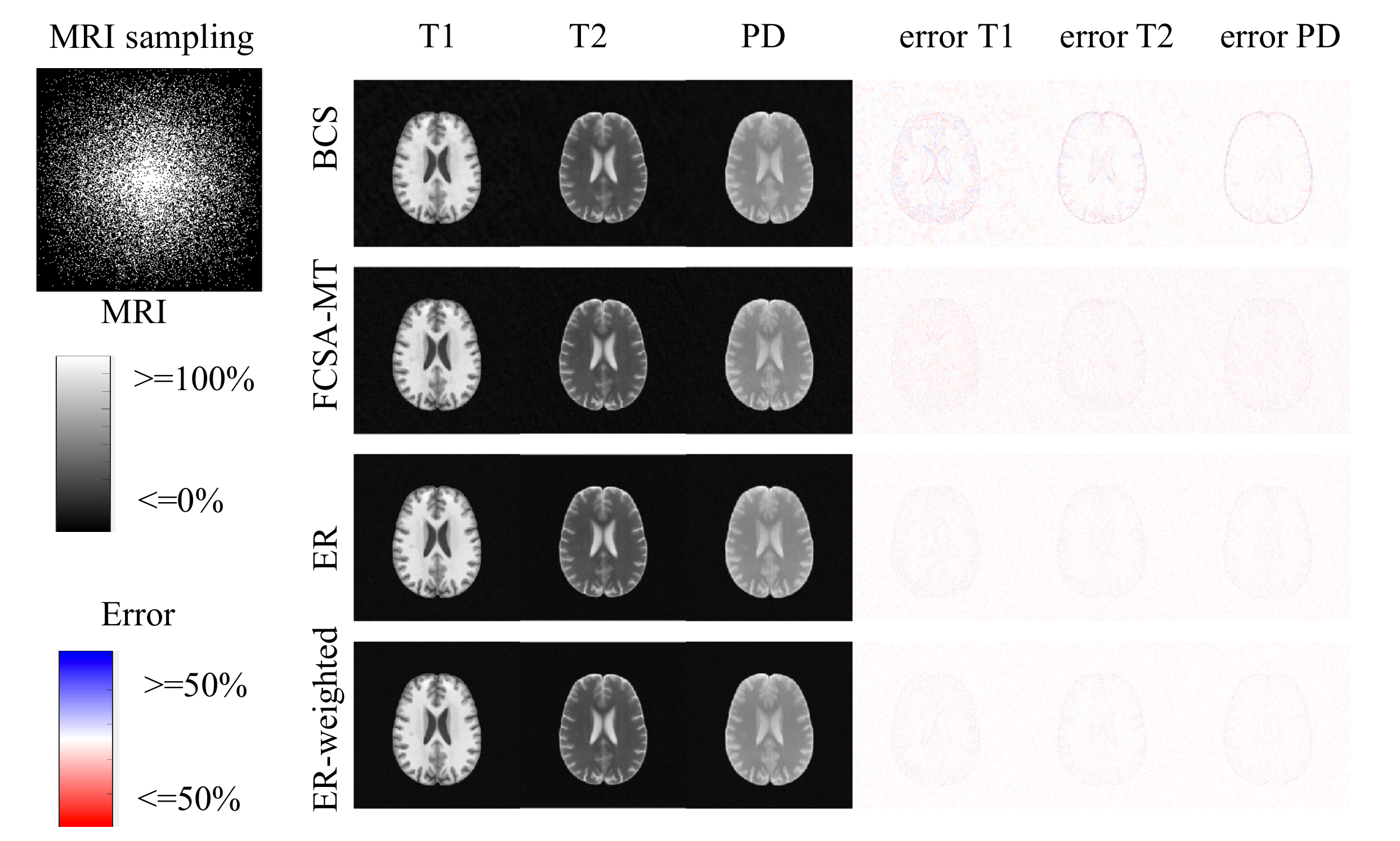}
\caption{Undersampling mask, reconstructed image, and error image for the BrainWeb image 
with noise level $\sigma=4$.} 
\label{fig:brainweb_sigma4_visual}
\end{figure}

We also conduct the same test using the BrainWeb image
with a Poisson mask of sampling ratio 25\%, also with noise level
$\sigma=4$ and $10$. The results are shown in Figure \ref{fig:brainweb_sigma_error_time}.
For this ``near real'' image, our method again shows significant improvement
of efficiency compared to other methods. In particular, ER-weighted
outperforms all methods in both cases, suggesting its superior 
robustness in multi-contrast image reconstruction.
The final reconstruction images and error images, similar to those for Shepp-Logan image,
are also plotted in Figure \ref{fig:brainweb_sigma4_visual} for the $\sigma=4$ case.

\begin{figure}[t]
\centering
\includegraphics[width=.32\textwidth]{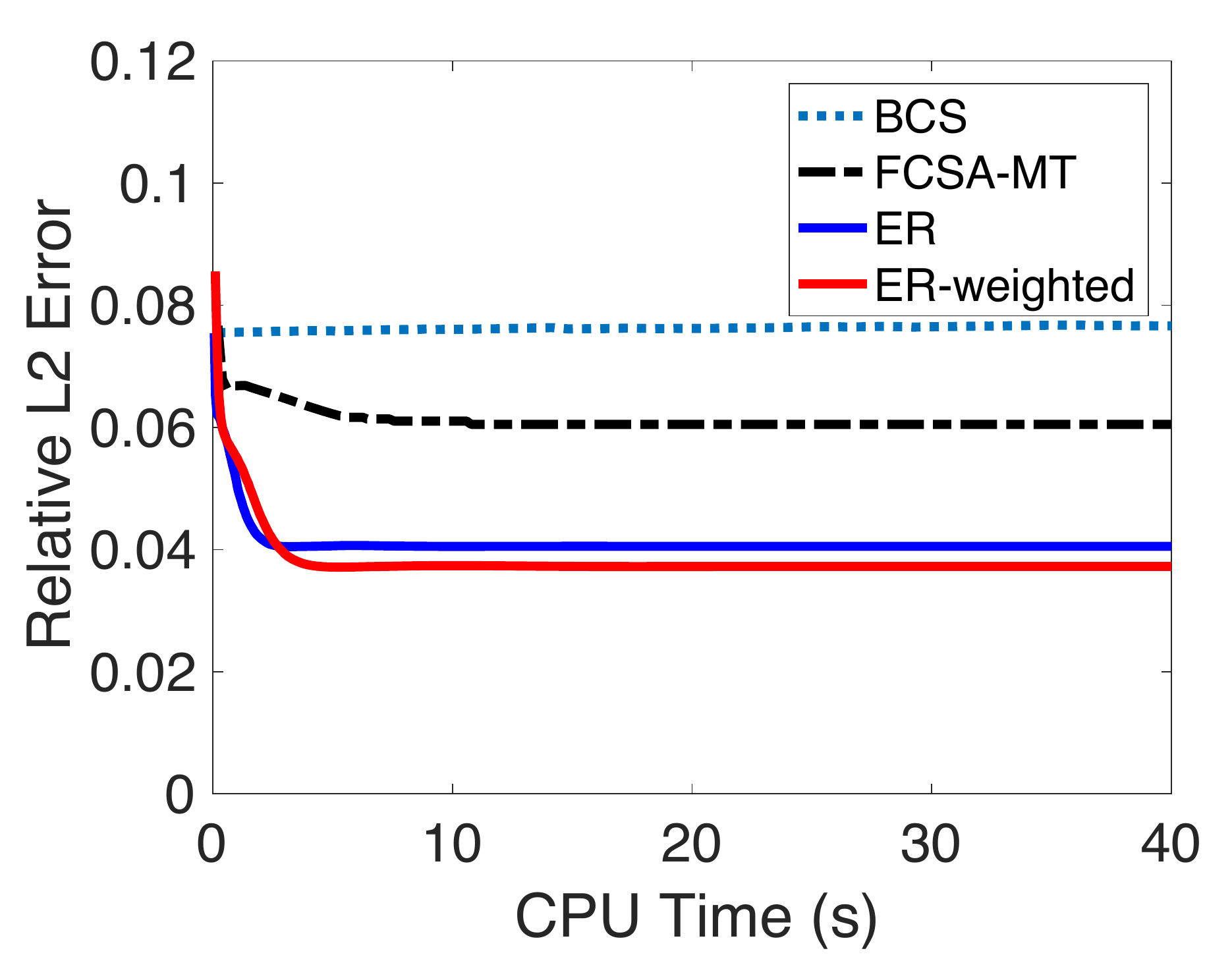}
\includegraphics[width=.32\textwidth]{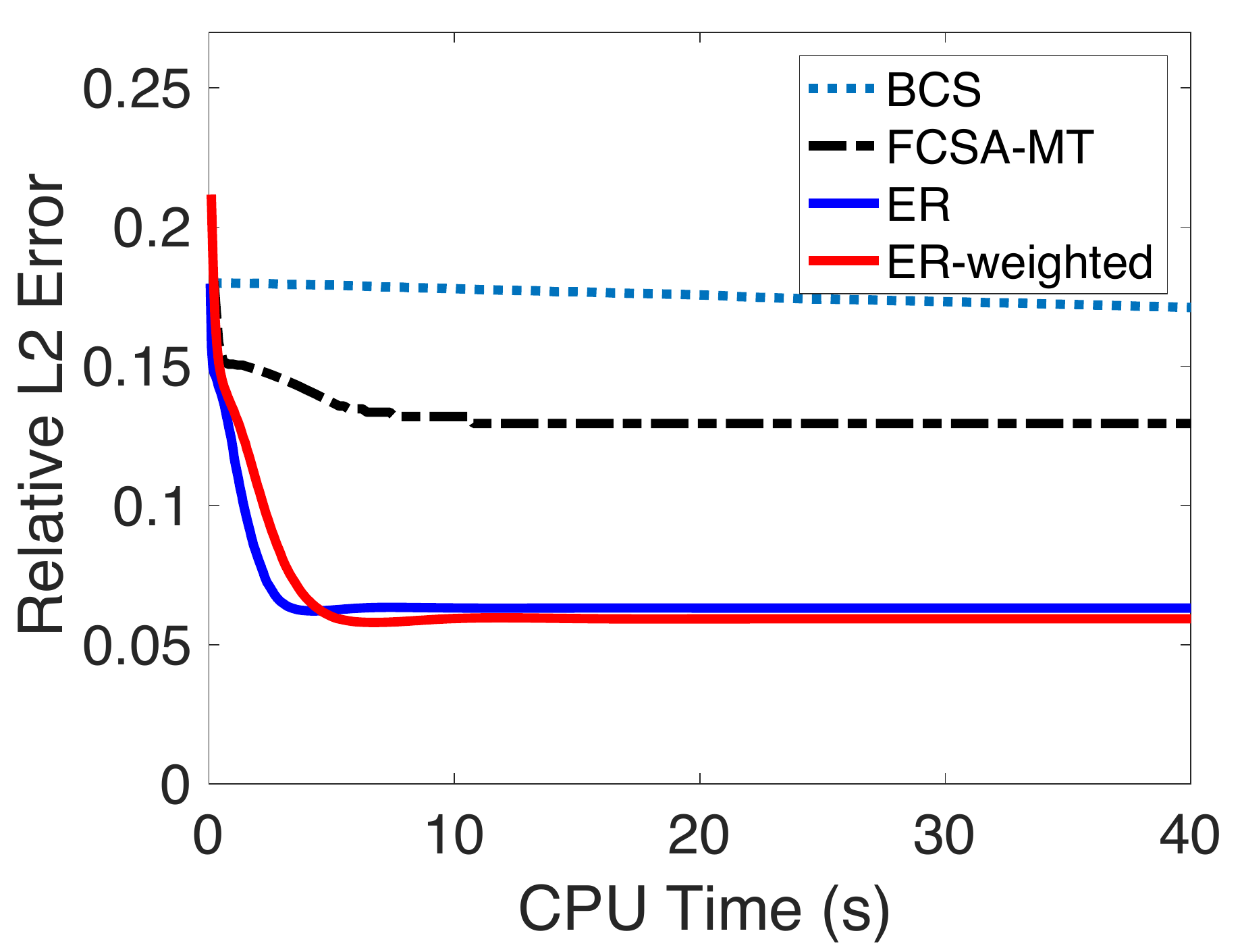}
\includegraphics[width=.32\textwidth]{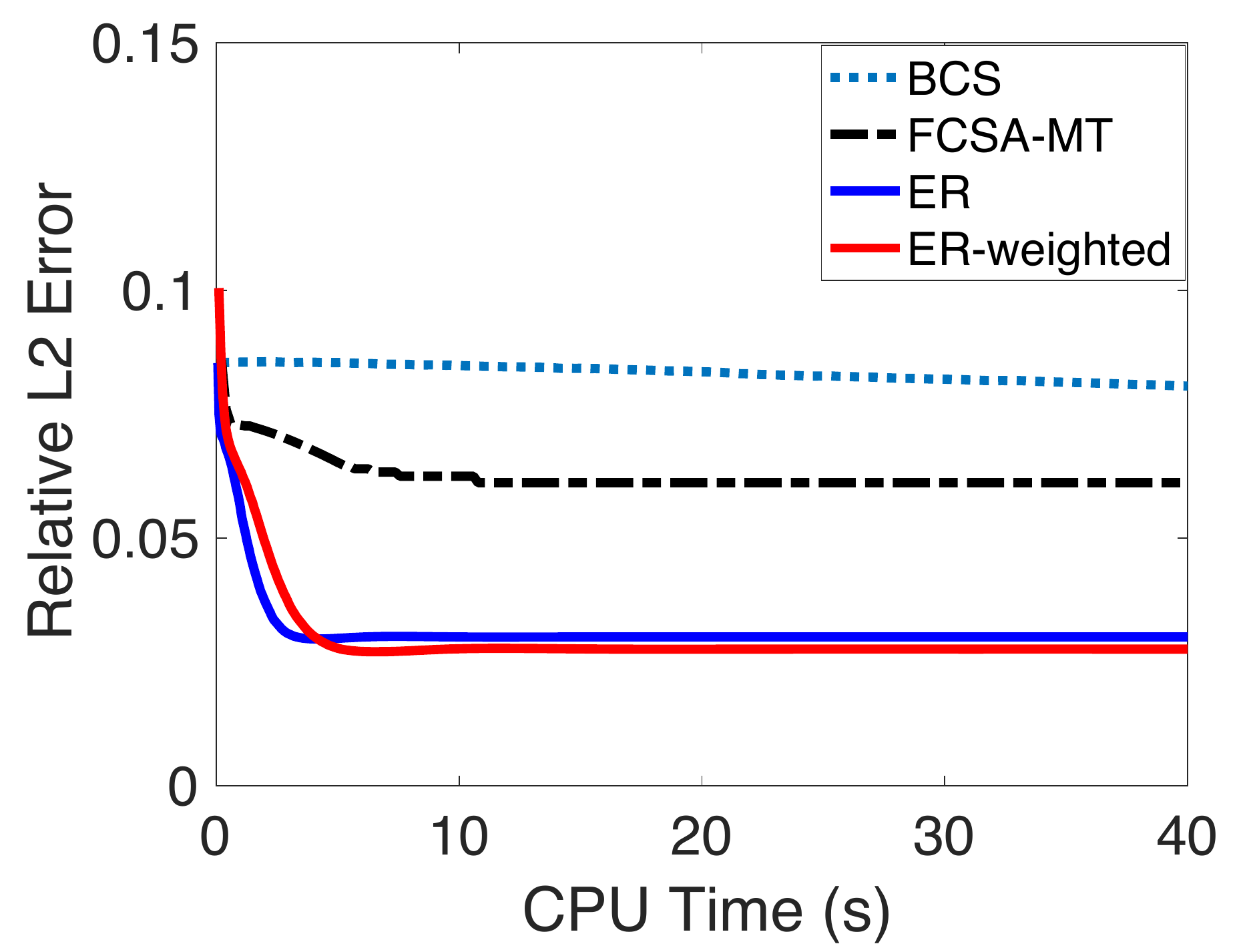}
\includegraphics[width=.32\textwidth]{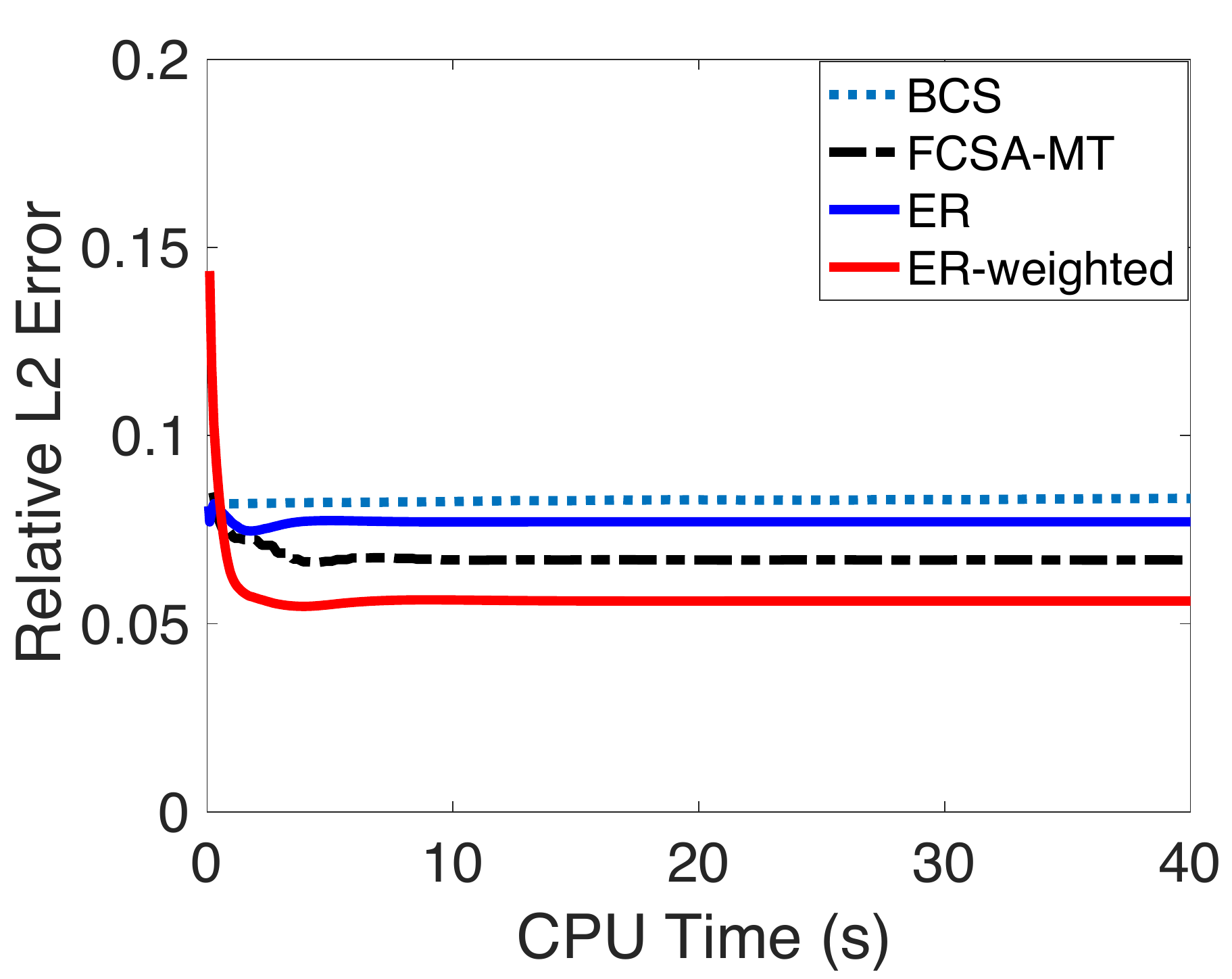}
\includegraphics[width=.32\textwidth]{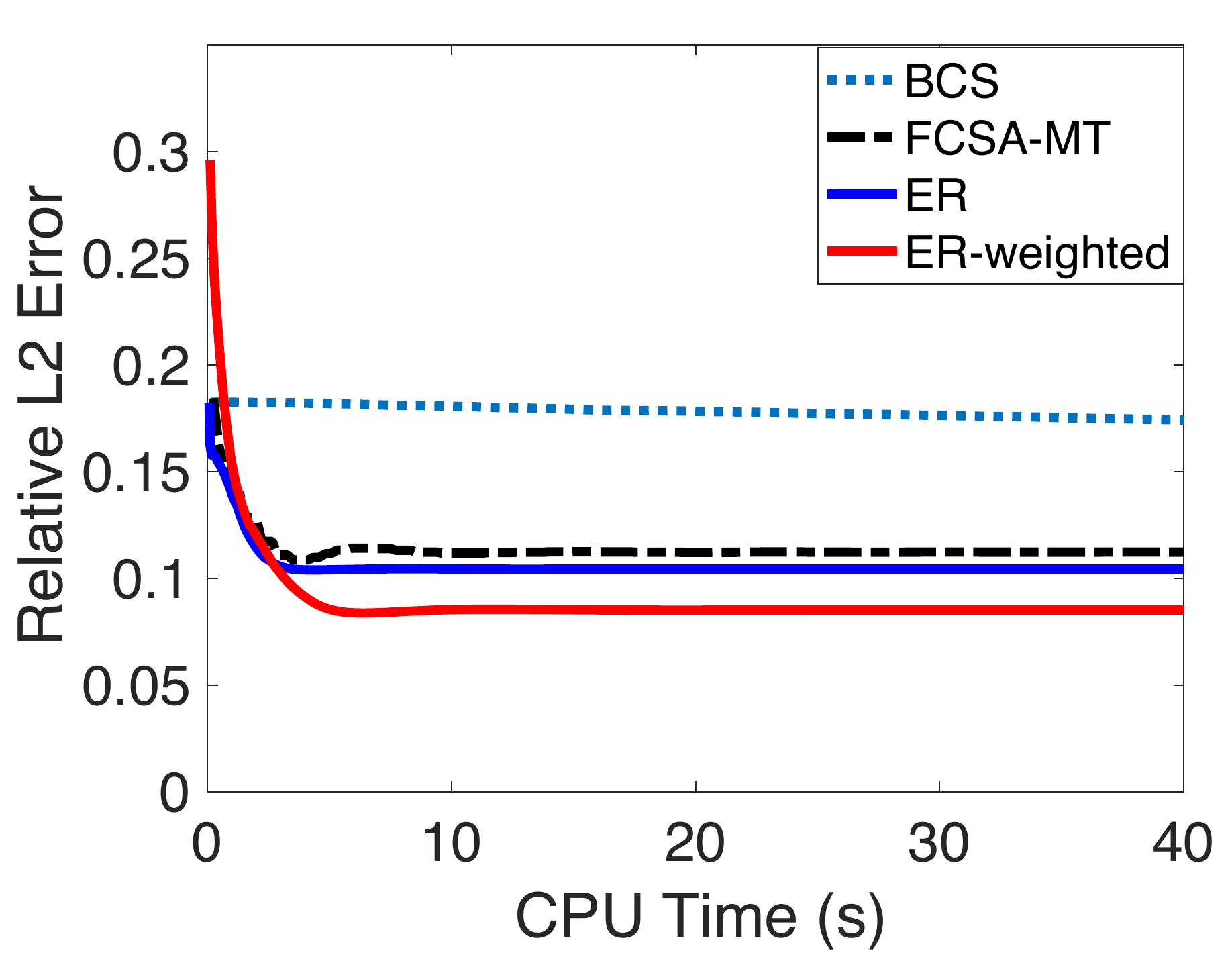}
\includegraphics[width=.32\textwidth]{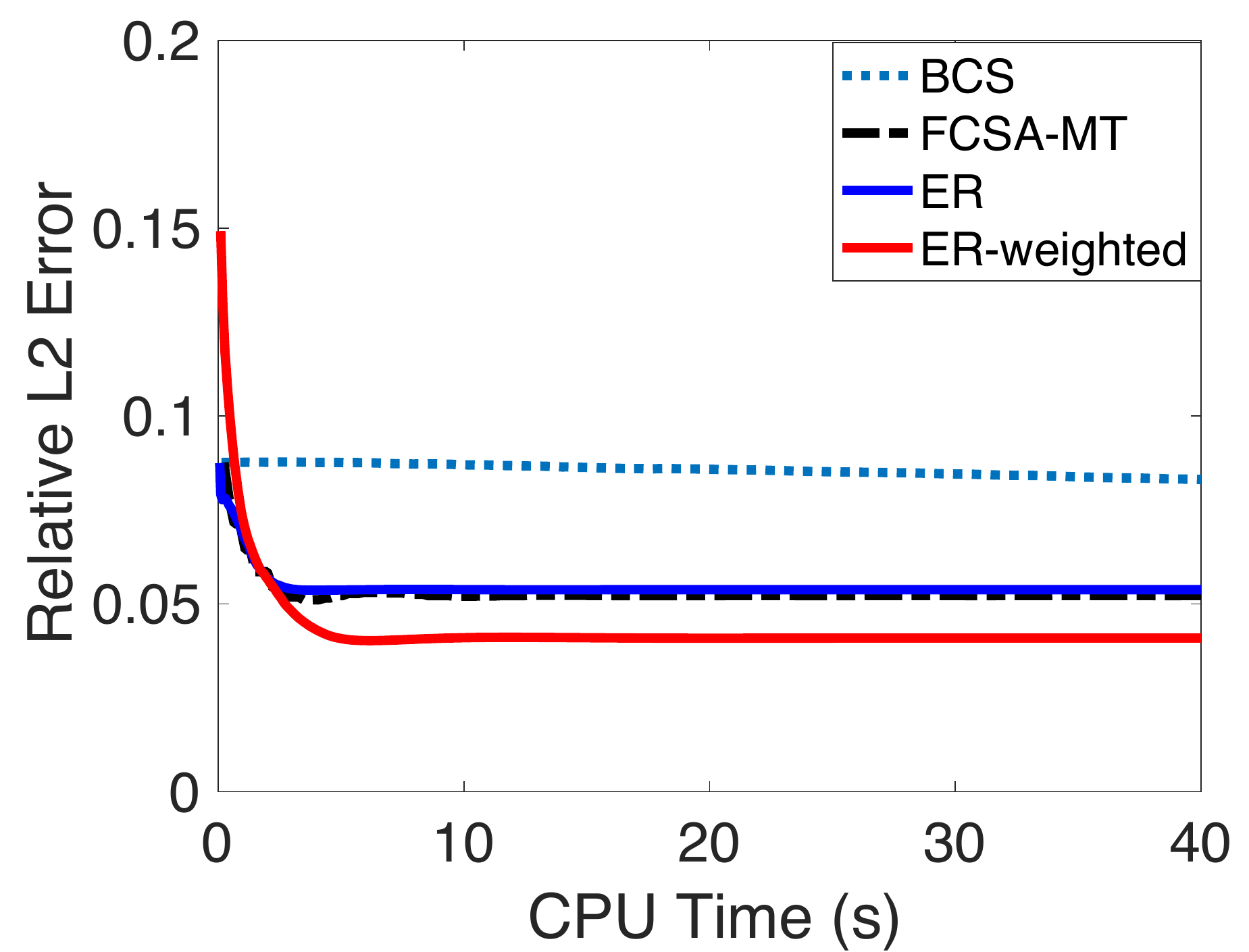}
\caption{Relative error vs CPU time of the comparison algorithms 
on in-vivo brain image with noise standard deviation  $\sigma=4$ (top row)
and $\sigma=10$ (bottom row).} 
\label{fig:brain_sigma_error_time}
\end{figure}

\begin{figure}[h!]
\centering
\includegraphics[width=.9\textwidth]{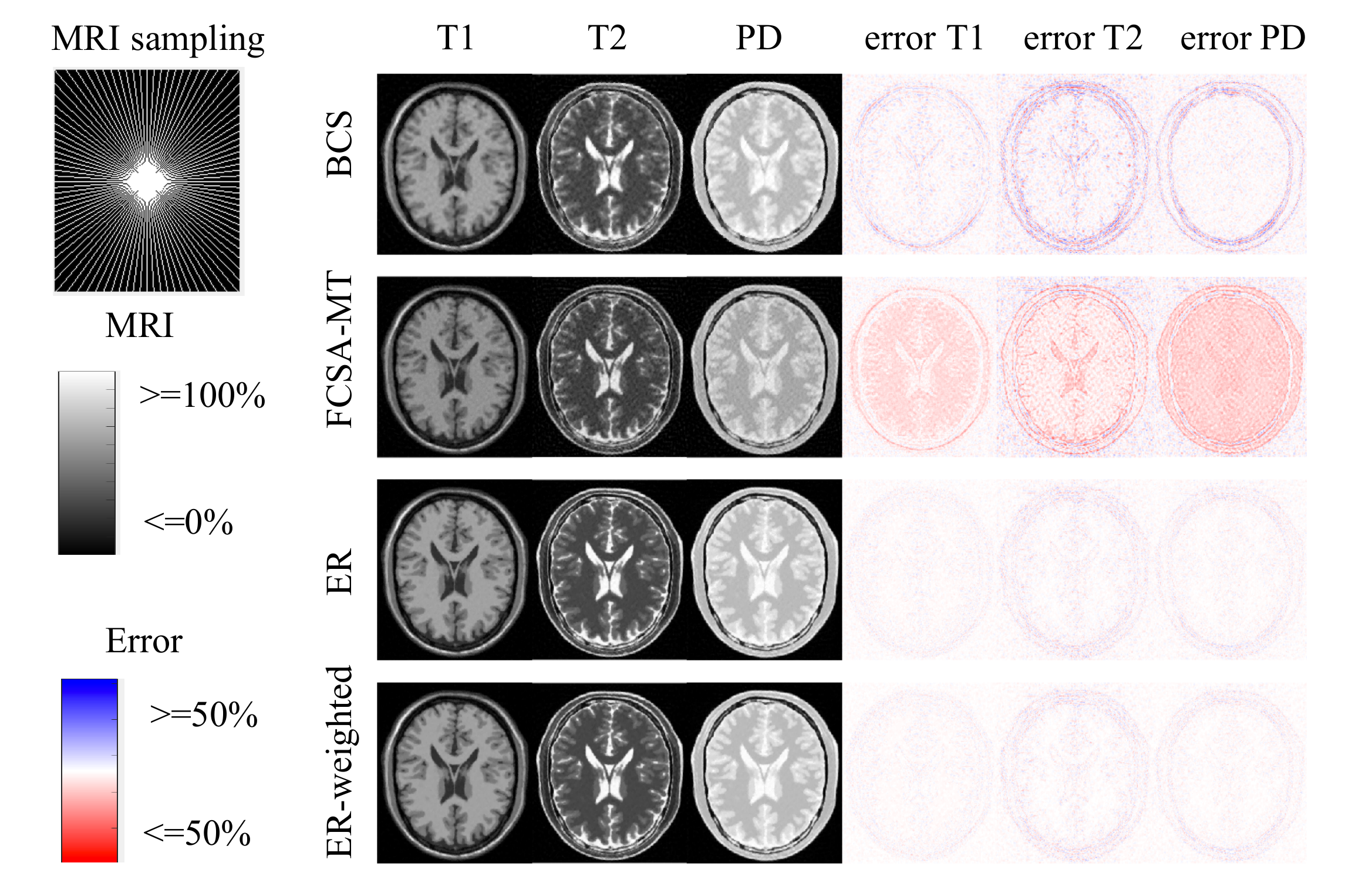}
\caption{Undersampling mask, reconstructed image, and error image for the in-vivo brain image 
with noise level $\sigma=4$.} 
\label{fig:brain_sigma4_visual}
\end{figure}

Finally, we conduct the same test on the in-vivo brain image
with radial mask of sampling ratio 25\%.
The relative error vs CPU time plots are shown in Figure 
\ref{fig:brain_sigma_error_time} ($\sigma=4$ and $10$).
For this image, our method further shows its promising performance
with significant improvement when compared to existing methods. 
The final reconstruction images and error images 
are also shown in Figure \ref{fig:brain_sigma4_visual} for the $\sigma=4$ case.
%
%

\section{Conclusion}
\label{sec:conclusion}
We proposed to reconstruct Jacobian of mutli-modal/contrast/channel image
from which we can resemble the underlying image.
We showed the relation between Jacobian and the observed data
when the underlying transform is Fourier, and formulate 
the reconstruction problem of Jacobian as an $l_1$ minimization.
Our new method then exhibits an optimal $O(1/k^2)$
convergence rate which outperforms the $O(1/k)$ rate of primal-dual based 
algorithms. We also derived closed-form solutions 
for the minimization subproblem as matrix-valued shrinkage.
The per-iteration complexity is thus very low.
Numerical results demonstrated the promising performance of the proposed method
when compared to the state-of-the-arts joint image reconstruction methods.

\section*{Acknowledgement}
The authors would like to thank Yao Xiao and Yun Liang for assisting Fang to run the MATLAB code and collect numerical results. This research is supported in part by National Science Foundation grant DMS-1719932 (Chen), IIS-1564892 (Fang), DMS-1620342 (Ye) and CMMI-1745382 (Ye), and National Key Research and Development Program of China No: 2016YFC1300302 (Fang) and National Natural Science Foundation of China No: 61525106 (Fang).

\appendix

\section*{Appendix}
\renewcommand{\thesubsection}{\Alph{subsection}}

\subsection{Computation of matrix-valued shrinkage}\label{apd:shrinkage}
We show that the minimization problem \eqref{eq:mtxshrink} has closed-form solutions
when the matrix $\star$-norm is Frobenius, induced 2-norm, or nuclear norm.

\begin{itemize}
\item \textbf{Frobenius norm.} 
In the case, all matrices can be considered
as vectors and hence the vector-valued shrinkage formula can be directly applied.
We omit the derivations here.

\item \textbf{Induced 2-norm.}
In the case, \eqref{eq:mtxshrink} becomes
\begin{equation}\label{eq:mtx2shrink}
\min_X \|X\|_2+\frac{1}{2\alpha}\|X-B\|_F^2 = \min_X\max_{\|\xi\|=\|\eta\|=1}\xi^T X \eta + \frac{1}{2\alpha}\|X-B\|_F^2
\end{equation}
where $\xi=(\xi_1,\xi_2)^T\in\Rbb^2$ and $\eta=(\eta_1,\dots,\eta_m)^T\in\Rbb^m$.

Switching min and max, and solving for $X$ with fixed $\xi$ and $\eta$, we obtain
$X=B-\alpha \xi \eta^T$ and the dual problem of \eqref{eq:mtx2shrink} becomes
\begin{equation}
\max_{\|\xi\|=\|\eta\|=1}\xi^T B \eta-\alpha \xi^T (\xi \eta^T) \eta + \frac{\alpha}{2}\|\xi\eta^T\|_F^2
=\max_{\|\xi\|=\|\eta\|=1}\xi^T B \eta- \frac{\alpha}{2}
\end{equation}
where we used the facts that $\xi^T (\xi \eta^T) \eta=\|\xi\|^2\|\eta\|^2=1$
and $\|\xi\eta^T\|_F^2=\sum_{i}\sum_{j}(\xi_i\eta_j)^2=\sum_{i}(\xi_i)^2\sum_{j}(\eta_j)^2=1$.
Therefore, $\xi$ and $\eta$ are the left and right singular vectors of $B$,
and hence the optimal solution of \eqref{eq:mtx2shrink} is $X^*=B-\alpha \xi\eta^T$.

To obtain $X^*$, one needs to compute the largest singular vectors of $B$, or the SVD of $B$.
Note that $B$ has size $2\times m$, and hence a reduced SVD for is sufficient.
In particular, for $m=2$, there is a closed form expression to compute SVD of $X$.
In general, SVD of a $2\times m$ matrix is easy to compute and has at most two singular
values.

\item \textbf{Nuclear norm}
This corresponds to standard nuclear norm shrinkage: compute SVD of $B=U\Sigma V^T$,
and set $X^*=U\max(\Sigma-\alpha,0)V^T$. Computation complexity is comparable
to the 2-norm case above.

\end{itemize}

\bibliographystyle{abbrv}
\bibliography{library}
\end{document}